\documentclass[12pt,reqno]{amsart}
\usepackage{amsmath}  
\usepackage{amsmath,amssymb}  
\usepackage[margin=1in]{geometry}
\usepackage{epstopdf}
\usepackage{graphicx}
\usepackage{tikz}
\usetikzlibrary{arrows}

\newtheorem{theorem}{Theorem} 
\newtheorem{proposition}{Proposition}  

\newtheorem*{corollary}{Corollary}
\newtheorem*{lemma*} {Lemma}

\def\R{\Bbb R} \def\Rp{\R_{>0}} \def\Z{\Bbb Z}  \def\N{\Bbb N}  \def\Q{\Bbb Q}  \def\Qp{\Q_{>0}} \def\CC{\Bbb C}
\def\={\;=\;}  \def\:{\;:=\;}   \def\+{\,+\,}  \def\m{\,-\,}  \def\p{\partial}
\def\iv{^{-1}} \def\h{\frac12} \def\sgn{\text{\rm sgn}}
\def\a{\alpha}  \def\b{\beta}  \def\l{\lambda}     \def\v{\varepsilon} \def\th{\vartheta} 
\def\G{\Gamma} \def\g{\gamma}   \def\s{\sigma} \def\D{\Delta}   \def\ep{\epsilon} 
 \newcommand{\sm}{\,\hbox{\tikz{\draw[line width=0.6pt,line cap=round] (3pt,0) -- (0,6pt);}}\,}  
\def\CD{\mathcal D}  \def\CR{\mathcal R}   \def\CW{\mathcal W}  \def\EE{\mathcal E}  \def\HH{\mathfrak H}
  
\def\sums{\sideset{}{^*}\sum}  \def\SL{\text{SL}(2,\Z)}
\def\sma#1#2#3#4{\bigl(\smallmatrix#1&#2\\#3&#4\endsmallmatrix\bigr)} 
 
\def\endpf{$\quad\square$}  
\def\wt{\widetilde} \def\wh{\widehat} \def\la{\langle}  \def\ra{\rangle} 
 
\def\be{\begin{equation}}  \def\ee{\end{equation}}  \def\bes{\begin{equation*}}  \def\ees{\end{equation*}}
\def\ba{\be\begin{aligned}}  \def\ea{\end{aligned}\ee}  \def\bas{\bes\begin{aligned}}  \def\eas{\end{aligned}\ees}

\def\r{\rho}  \def\sg{\sigma}  \def\C{{\sf C}}  \def\hC{\wh\C}  
\def\GR{\sf{GRH}}     
\def\re{\mathrm{Re}\,}  
\def\H{G}    
\def\t{\tau} \def\ch{\chi}    \def\V{\mathcal V}  \def\wV{\wh{\V}}
\def\dst{{\sf d}}  \def\mon{{\sf p}}      \def\Hd{\mathrm{H}}
      \def\A{{\sf A}} 
\def\L{L}   

\begin{document}

\title[Cotangent sums and the generalized Riemann hypothesis]
 {Cotangent sums, quantum modular forms, \\ and the generalized Riemann hypothesis}
\author{John Lewis}  \author{Don Zagier}    \maketitle

\section{Introduction and main results} \label{sec.intro}

To answer immediately a question that the reader may be asking, we should say from the outset that the
generalized Riemann hypothesis (henceforth, just GRH) of the title refers only to the $L$-series
associated to odd real Dirichlet characters and that we do not aim to prove it, but merely
to show its equivalence to an asymptotic statement about the determinants of certain matrices
whose entries $c_{m,n}$ are given by finite sums of cotangents of rational multiples of~$\pi$.  To answer
a second natural question, since it may seem unusual for an author to contribute to the proceedings volume
of a conference partly intended to celebrate his own birthday, we should also say a few words about the genesis 
of this joint paper. Both the definition of the cotangent sums considered and their relation to~GRH were found 
by the first author and presented in his talk at the conference in question.  In the course of that talk he 
mentioned the homogeneity property $c_{\ell m,\ell n}=\ell c_{m,n}$, leading the second author to ask whether the 
function $C:\Q\to\R$ defined by $C(m/n)=c_{m,n}/n$ might be related to a quantum modular form in the sense 
of~\cite{QMF}. The answer to this question turned out to be positive.  Furthermore, and unexpectedly, this quantum 
modular aspect turned out also to be related to the earlier joint work of the authors~\cite{LZ} on the analogues 
of classical period polynomials for Maass wave forms. A glance at the title of the final section of this paper will 
show why it seemed natural, and even irresistible, to publish both parts of the story in these proceedings.  

We now describe the main results of the paper in more detail. We denote by
\be\label{chi4L4}
 \chi(n)\,=\,\begin{cases} (-1)^{(n-1)/2} & \text{if $n$ is odd,} \\ 0 & \text{if $n$ is even,} \end{cases}
 \qquad L(s)\,=\, 1 \m \frac1{3^s} + \frac1{5^s} \m \frac1{7^s} \+\cdots\,. \ee
the primitive Dirichlet character of conductor~4 and its associated $L$-series. (In Section~\ref{sec:OtherChi} 
we will describe the generalization to other odd primitive real Dirichlet characters.)  
For positive integers $m$ and~$n$ we define a real number $c_{m,n}$ by 
\be\label{defcmn}
 c_{m,n} \: \frac4\pi\,\sum_{j,\,k>0} \frac{\chi(j)\,\chi(k)}{\max(j/m,k/n)}\,, 
\ee
where the conditionally convergent sum is to be interpreted as the limit for $N\to\infty$ of the sum 
over the rectangle $1\le j\le mN$, $1\le k\le nN$. These numbers are studied in more detail in~\S\ref{Sec2}, 
where we will see that they are algebraic numbers (in fact, algebraic integers) that can be expressed as 
finite sums of cotangents of odd multiples of $\pi/4m$ or~$\pi/4n$, a typical example being
$$ c_{3,4} \= \cot\Bigl(\frac\pi{12}\Bigr)  \+ \cot\Bigl(\frac{5\pi}{12}\Bigr) 
  \+ \cot\Bigl(\frac{9\pi}{16}\Bigr) + \cot\Bigl(\frac{13\pi}{16}\Bigr) \= 6 \m 2\,\sqrt{2+\sqrt2}\;. $$
For each $N\in\N$ we define the two symmetric matrices 
\be\label{defCN}  \C_N \: \begin{pmatrix} c_{1,1} & \cdots & c_{1,N} \\ \vdots & \ddots & \vdots \\ 
  c_{N,1} & \cdots & c_{N,N} \end{pmatrix}\;, \qquad 
\hC_N \: \left(\begin{array}{ccc|c}   &  & & 1 \\ [-4pt] & \C_N & & \vdots   \\  & & & N \\
       \hline   1 & \cdots & N & 4N/\pi \\[1pt]  \end{array}\right)\,, \ee
which we will show later are positive definite, and define a function $R:\N\to\R$ by
\be \label{defr}  R(N) \: N\,\frac{\det \C_N}{\det\hC_N}\;. \ee
The first main result of this paper is then the following.
\begin{theorem} \label{GRH}
The following two statements are equivalent.
\begin{itemize} \item[(a)] The function $R(N)$ is unbounded.
\item[(b)] The $L$-series $L(s)$ has no zeros in the half-plane $\sg:=\Re(s)>\frac12$~(GRH).
\end{itemize}  \end{theorem}
\noindent The proof of this theorem will be given in \S\ref{sec:GRH}. The direction (a)$\,\Rightarrow\,$(b) is fairly 
elementary, while the reverse direction uses functional analysis and techniques coming from old work of Beurling. 
The result that we actually prove (Theorem~\ref{Thm3} in~\S3) is in fact stronger than Theorem~\ref{GRH} in both 
directions: if GRH holds, then $R(N)$ is not merely unbounded, but actually tends to infinity as $N\to\infty$, and if 
$L(s)$ has a zero $\r$ with $\Re(\r)>\frac12$, then the function $R(N)$ is bounded by a number depending explicitly on~$\r$. 
Of course we believe that the two statements in Theorem~\ref{GRH} are both true rather than both false, the evidence being
on the one hand the generally held belief in the validity of~(b) and on the other hand the following plot of the function 
$R(N)$ for $N\le100000$, which suggests that $R(N)$ tends to infinity and perhaps even grows roughly linearly with~$\log N$.  
(The straight-line fit shown in the picture is the graph of the function $\,5.18\,\log N\,$.)

\medskip
\begin{center}
  \includegraphics[angle=0,height=8cm]{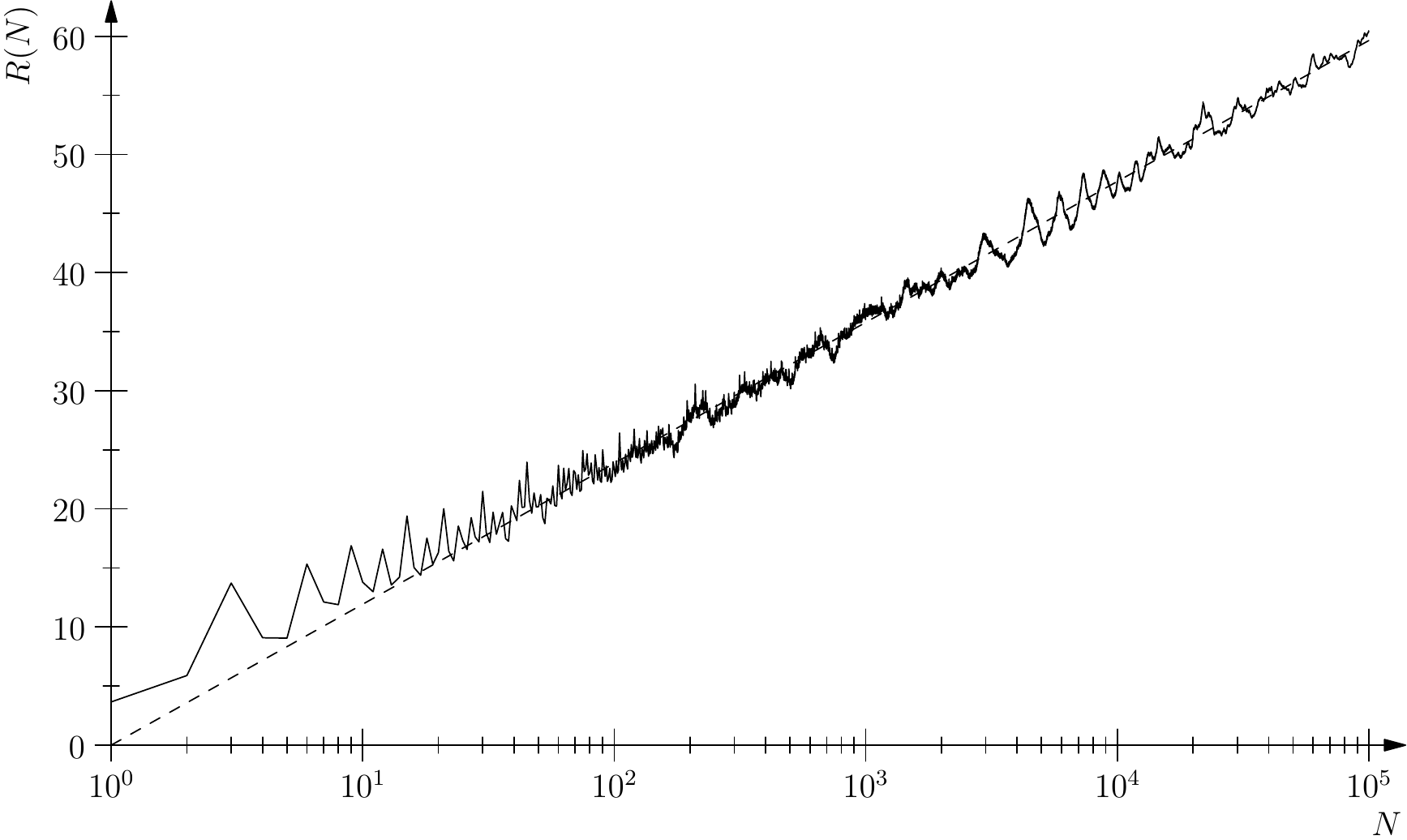} \\
  \textbf{Figure 1. Graph of the function $R(N)$}
\end{center}

\medskip

The second main result concerns the quantum modular nature of certain functions related to the 
coefficients~$c_{m,n}$.  On the one hand, the numbers $c_{m,n}$ defined by~\eqref{defcmn} have 
the obvious homogeneity property $c_{\ell m,\ell n}=\ell\,c_{m,n}$ and can be rewritten in the form
\be \label{defC}  c_{m,n}\= \frac{4n}\pi\,\,C\Bigl(\frac mn\Bigr)\,, \qquad
 C(x) \: \sum_{j,\,k>0} \chi(j)\,\chi(k)\,\min\Bigl(\frac xj,\,\frac1k\Bigr)\qquad(x\in\Q_{>0})\,. 
\ee
On the other hand, by splitting the sums in~\eqref{defcmn} or~\eqref{defC} into two pieces depending 
on which of its two arguments realizes the ``max" or ``min", we obtain the decompositions
\be\label{defhmn} c_{m,n} \= h_{m,n}\+h_{n,m}\,, \qquad
 C(x) = H(x)\+x\,H\Bigl(\frac1x\Bigr)
\ee
where the function $H:\Q_{>0}\to\R$ and numbers $h_{m,n}\in\R$ are defined by
\be \label{defH}  
H(x) \: \sums_{0<j\le kx}\frac{\chi(k)\,\chi(j)}k\,,\qquad
 h_{m,n} \:  \frac{4n}\pi\, H\Bigl(\frac mn\Bigr)\,.  \ee
(Here and later an asterisk on a summation sign means that terms with equality---in this case, those with 
$j=kx$---are to be counted with multiplicity 1/2.) We will show in Section~\ref{Sec2} that $h_{m,n}$ (or its sum 
with~1/2 if $mn$ is odd) is a sum of cotangents of odd multiples of~$\pi/4n$ and hence is an 
algebraic integer, refining the corresponding properties for $c_{m,n}$ given above.
 
We can now come to the quantum modular form property mentioned in the
opening paragraph.  The graph of the function~$H$ on~$\Q$ looks as follows:

\bigskip

\begin{center}
  \includegraphics[width=0.8\linewidth]{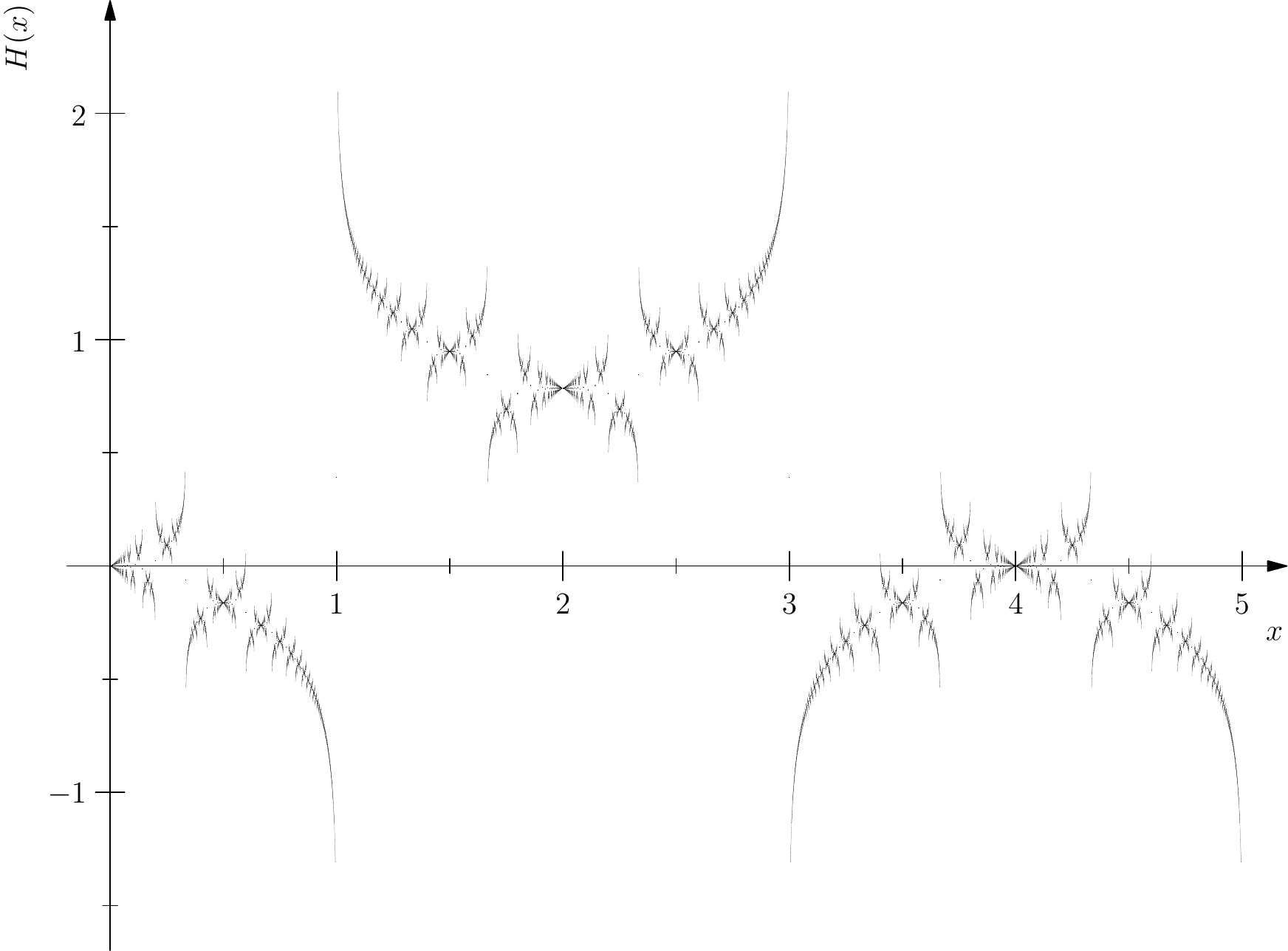}\\
  \vspace{2mm}
  \textbf{Figure 2. Graph of the function $H(x)$}
\end{center}

\bigskip \noindent
This graph suggests that the function $H$ is only naturally defined on~$\Q$ and is not well-behaved as
a function on~$\R$, having a possibly dense set of discontinuities and certainly not differentiable at generic 
points. The graph of the function $C$ shown in Figure~3, on the other hand, suggest that it is 
everywhere continuous and possibly even differentiable except at rational points with odd numerator and 
denominator, where it seems to be left and right differentiable but with a discontinuous derivative.  
Our second main result, proved in Section~\ref{sec:anal}, refutes the last of these statements and confirms the others.

\medskip
\begin{center}
  \includegraphics[width=0.8\linewidth]{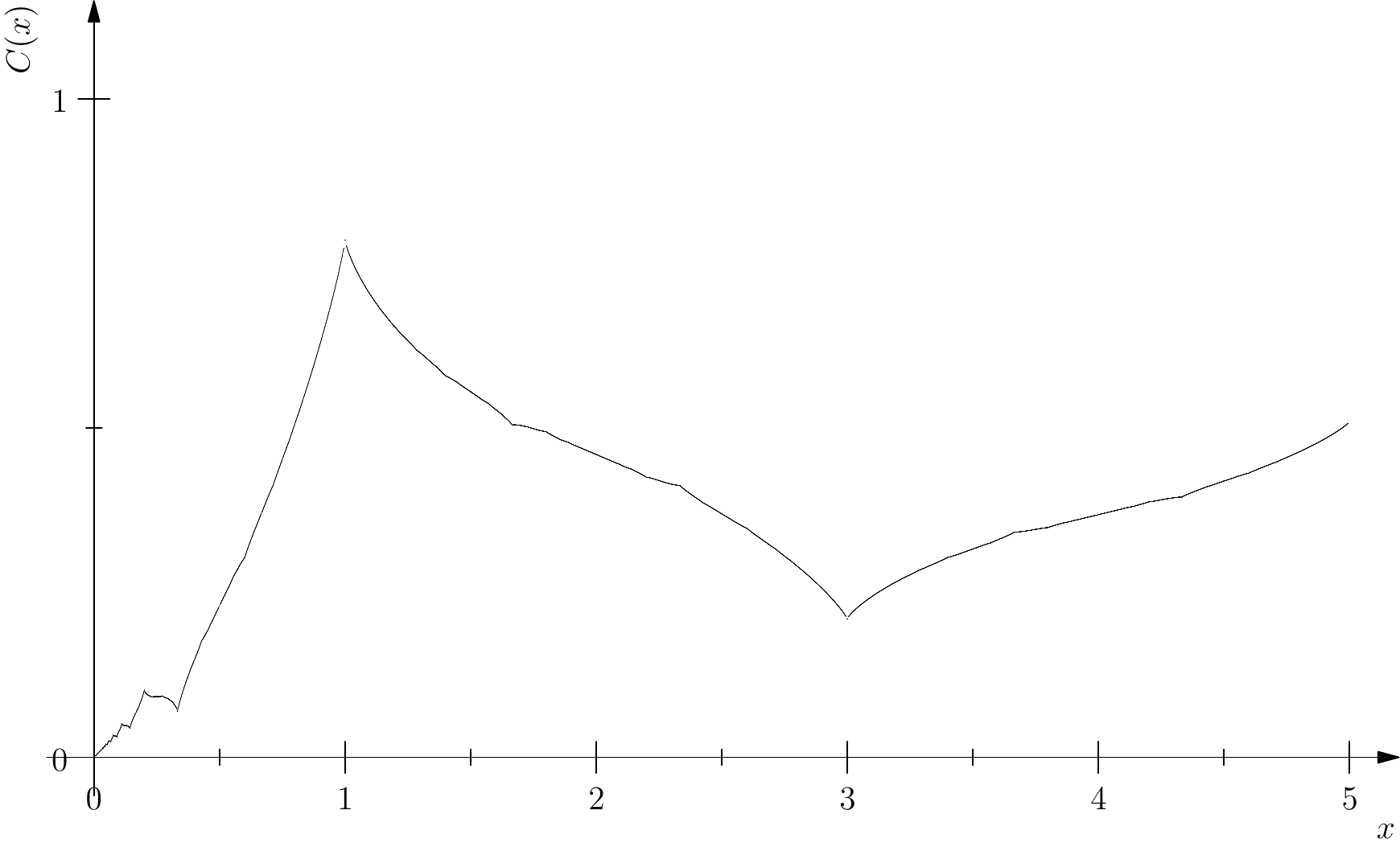}\\
  \vspace{0mm}
  \textbf{Figure 3. Graph of the function $C(x)$}
\end{center}

\bigskip

\begin{theorem} \label{continuity}
The function $H(x)$ is discontinuous at all rational points~$x$ with odd numerator and denominator, 
and is continuous but not left or right differentiable at other rational points. The function $C(x)$ is 
continuous everywhere, is not left or right differentiable at rational points with odd 
numerator and denominator, but is differentiable at all other rational points. 
\end{theorem} 

The similarity between the continuity and differentiability statements for the two functions 
$H(x)$ and~$C(x)$ in this theorem is not coincidental, since we will see in Section~\ref{sec:anal}
that the derivative of $C(x)$ is equal to $H(1/x)$.

We now discuss the ``quantum modular form" aspect of Theorem~\ref{continuity}.
Note first that the function~$H$ is periodic of period~4 and hence extends to a periodic 
function on all of~$\Q$. This function is even, so we could also extend $H$ from $\Qp$ to~$\Q$
by $H(x):=H(|x|)$. One might then think of extending $C$ to negative values of the arguments
by setting $C(x)=H(x)+xH(1/x)$ for all $x\ne0$, but if we did that then $C(x)$ for $x<0$ would
exhibit oscillatory behavior of the same sort as~$H(x)$. Instead, we extend $C$ to all of~$\Q$
as an even function, i.e., we set $C(x)=C(|x|)$ for $x\ne0$ and extend this by continuity 
to~$C(0)=0$.  Then the relationship between $C$ and $H$ becomes
\be\label{CH} H(x) \+ |x|\,H(1/x) \= C(x) \qquad(x\ne0)\,. \ee
On the other hand, the 4-periodicity property of $H(x)$ can be strengthened to the statement
\be\label{Hper2}  H(x)\+ H(x+2)\= \frac\pi4\,, \ee
saying that up to sign and up to a constant $H$ actually has period 2. Also, since $H$ is 
even, we can replace the term $H(1/x)$ in~\eqref{CH} by $H(-1/x)$. But the two matrices
$S=\sma0{-1}10$ and $T^2=\sma1201$ generate the subgroup~$\G_\th$ of index~3 in the full modular
group $\,\SL=\langle S,T\rangle$ consisting of matrices congruent to $\sma1001$ 
or $\sma0110$ modulo~2 (this is the so-called ``theta group," under which the Jacobi theta 
function $\th(z)=\sum_ne^{\pi in^2z}$ transforms like a modular form of weight~1/2), so
equations~\eqref{CH} and~\eqref{Hper2} can be combined to the following statement: 
\begin{corollary} \label{QMF}
The function $H:\Q\to\R$ satisfies the transformation property
$$  \v(\g)\,|cx+d|\,H\Bigl(\frac{ax+b}{cx+d}\Bigr) \= \,H(x) \;+\;\text{\rm(continuous function)} $$
for every matrix $\g=\sma abcd\in\G_\th$, where $\v:\G_\th\to\{\pm1\}$ is the homomorphism sending
both generators $\sma1201$ and $\sma0{-1}10$ of~$\G_\th$ to~$-1$. 
\end{corollary}

This corollary tells us precisely that the function $H$ is a quantum modular form in the sense of~\cite{QMF} 
(except that there is an absolute value sign here that was not present there). Recall that a ``quantum 
modular form" is by definition a function on~$\Q$ that does not extend nicely to~$\R$ but which modulo (the 
restriction from $\R$ to~$\Q$ of) better-behaved (here, continuous) functions on~$\R$ transforms like a modular 
form with respect to some subgroup of finite index of~$\SL$ (here the theta group~$\G_\th$) with some 
weight (here~$-1$) and character~(here~$\v$). In Section~\ref{sec:MF} of this paper we will show that the quantum 
modularity property of $H$ is related to an actual modular form (specifically, a Maass Eisenstein
series of eigenvalue~$1/4$) on the group~$\G_\th$.

The authors would like to thank Steven Johnson~(MIT) and Alexander Weisse~(MPI) for invaluable 
help with the numerical aspects of this paper.


\section{Properties of the numbers $c_{m,n}$ and $h_{m,n}$} \label{Sec2}

In this section we give an integral representation for the numbers $c_{m,n}$ that will be used
later both to extend the function~$C(x)$ continuously from~$\Q$ to~$\R$ and as the key to the proof 
in~Section~\ref{sec:GRH} of the relation between the numbers $c_{m,n}$ and the zeros of the Dirichlet 
$L$-series~$L(s)$.  We will also prove the formulas for $c_{m,n}$ and $h_{m,n}$ as algebraic numbers 
(finite sums of cotangents).

The first step is to rewrite the definition of~$H$ given in~\eqref{defH} more compactly in the form
\be\label{HS}  H(x) \= \sum_{k=1}^\infty \frac{\ch(k)}k\,S(kx) \qquad(x\in\Q)\,, \ee
where $S(x)$ is the function defined on $\Rp$ by
\be\label{defS} S(x) \: \sums_{0<j\le x}\ch(j)\,. \ee
The function~$S$, like $H$ itself, is periodic of period~4, and we can use this property to extend~it 
to all of~$\R$. (Equivalently, we can extend it by $S(-x)=S(x)$ and $S(0)=0$, 
or by observing that the summation condition in~\eqref{defS} can be replaced by $J<j\le x$
for any integer~$J<x$ divisible by~4, which works for all real~$x$.) The function~$S$ 
is then the step function given by
\be\label{step}  S(x) \= \begin{cases} 1 & \text{if $4k+1<x<4k+3$ for some $k\in\Z$,} \\
  0 & \text{if $4k-1<x<4k+1$ for some $k\in\Z$,} \\  \frac12 & \text{if $x$ is an odd integer,} \end{cases} \ee
as pictured in Figure~4.  Note that the sum in~\eqref{HS} converges (conditionally) for $x\in\Q$ by summation
by parts, because the function $\ch(k)S(kx)$ is then an odd periodic function of~$k$ and hence has 
average value zero, so that the sums $\sum_{k=1}^K\ch(k)S(kx)$ are bounded.

\begin{center}
  \begin{tikzpicture}[scale=1.8]
    \newcommand{\charfour}{
      \path[fill=black!10!white] (1,0) rectangle (3,1);
      \draw[very thick] (0,0) -- (1,0);
      \draw[dashed] (1,0) -- (1,1);
      \draw[very thick] (1,1) -- (3,1);
      \draw[dashed] (3,1) -- (3,0);
      \draw[very thick] (3,0) -- (4,0);
      \path[draw,fill=white] (1,0) circle [radius=1pt];
      \path[draw,fill=white] (3,0) circle [radius=1pt];
      \path[draw,fill=white] (1,1) circle [radius=1pt];
      \path[draw,fill=white] (3,1) circle [radius=1pt];
      \path[draw,fill=black] (1,0.5) circle [radius=1pt];
      \path[draw,fill=black] (3,0.5) circle [radius=1pt];
    }
    \begin{scope}
      \clip(-1.5,-0.5) rectangle (5.5,1.5);
      \begin{scope}[xshift=-4cm]
        \charfour
      \end{scope}
      \begin{scope}[xshift=0cm]
        \charfour
      \end{scope}
      \begin{scope}[xshift=4cm]
        \charfour
      \end{scope}
    \end{scope}
    \draw[->] (-1.5,0) -- (5.6,0);
    \draw[->] (0,-0.1) -- (0,1.3);
    \foreach \x in {-1,0,1,...,5} {
      \draw (\x,0.05) -- (\x,-0.05) node[anchor=north] {$\x$};
    }
    \draw (0.05,0.5) -- (-0.05,0.5) node[anchor=east] {$\tfrac{1}{2}$};
    \draw (0.05,1) -- (-0.05,1) node[anchor=east] {$1$};
  \end{tikzpicture}\\[2mm]
  \textbf{Figure 4. Graph of the function $S(x)$}
\end{center}

\bigskip\medskip

We use the bounded periodic function $S$ to prove the following formula for the $c_{m,n}\,$: 

\begin{proposition} \label{prop:cmnintegral} 
The numbers $c_{m,n}$ defined by~\eqref{defcmn} have the integral representation
\be\label{cmnintegral}  c_{m,n} \= \frac4\pi\,\int_0^\infty\!S(mt)\,S(nt)\,\frac{dt}{t^2} \qquad(m,\,n\in\N). \ee
\end{proposition}
\noindent{\it Proof.} From the definition of $S$ we find
$$\int_0^N\!S(mt)\,S(nt)\,\frac{dt}{t^2} \= \int_0^N\!\sums_{0<j\le mt\atop 0<k\le nt}\ch(j)\ch(k)\,\frac{dt}{t^2} 
\= \sum_{0<j\le mN\atop 0<k\le nN}\biggl(\frac{\ch(j)\,\ch(k)}{\max(j/m,k/n)} \m \frac{\ch(j)\,\ch(k)}N\biggr)\,,$$
and since the second sum is identically zero for $N$ divisible by~4, equation~\eqref{cmnintegral} follows by letting
$N$ tend to infinity in~$4\Z$.  \endpf

Our next result expresses $c_{m,n}$ as a finite sum of algebraic cotangent values indexed by lattice 
points in the rectangle $[0,m/2]\times[0,n/2]$ satisfying a certain ``selection rule."

\begin{proposition} \label{prop:selection}
Let $m$ and $n$ be positive integers.  Then the number defined by~\eqref{defcmn} is 
given as the following integral linear combination of cotangents of rational multiples of~$\pi$:
$$  c_{m,n} \,=  \sum_{ 0\le j\le m/2\atop 0\le k\le n/2} 
  \begin{cases} \cot(\pi a)-\cot(\pi b) & 
  \text{if $a:=\max\bigl(\frac{4j+1}{4m},\frac{4k+1}{4n}\bigr)<
  b:=\min\bigl(\frac{4j+3}{4m},\frac{4k+3}{4n},\frac12\bigr)\,,$} \\ 0 &\text{if $a\ge b$\,.}\end{cases}$$
\end{proposition}
\begin{corollary}  The numbers $c_{m,n}$ are algebraic integers for all $m,\,n>0$.
\end{corollary} 

Using this proposition we can calculate any $c_{m,n}$ in closed form as an algebraic number. 
For example, the first $4\times4$ matrix of values (cf.~equation~\eqref{defCN}) is given by
$$ \C_4 \= \begin{pmatrix} 1&2-\sqrt2&2-\sqrt3&2-2\sqrt{2-\sqrt2}\; \\
2-\sqrt2&2&\sqrt2&4-2\sqrt2\\2-\sqrt3&\sqrt2&3&6-2\sqrt{2+\sqrt2}\\
\,2-2\sqrt{2-\sqrt2}&4-2\sqrt2&6-2\sqrt{2+\sqrt2}&4\end{pmatrix}\;.$$

\smallskip
\noindent{\it Proof.} Since $S(x)$ is even and periodic of period~4, we can use Euler's formula
$\sum_{k\in\Z}(k+t)^{-2}=\pi^2/\sin^2(\pi t)$ to rewrite the integral representation~\eqref{cmnintegral} as
$$ c_{m,n} \=  \pi\,\int_0^{1/2} S(4mx)\,S(4nx)\,\frac{dx}{\sin^2(\pi x)}\,.$$
But from the definition of $S(4mx)$ and $S(4nx)$ as characteristic functions, it follows that the restriction
of $S(4mx)S(4nx)$ to~$[0,\frac12]$ is the characteristic function of the union of all intervals $[a,b]$
with $a<b$ as given in the proposition. This proves the proposition because 
$\pi\int_a^b\csc^2(\pi x)\,dx = \cot(\pi a)-\cot(\pi b)$.
The corollary then follows because by elementary algebraic number theory the number $\cot(\pi\a)$ is an 
algebraic integer for every rational number $\a$ having an even denominator. 
\endpf

We next give the corresponding formula for the numbers $h_{m,n}$, which is in fact simpler.
\begin{proposition} \label{prop:hmnformula}
Let $m$ and $n$ be positive integers.  Then the number defined by~\eqref{defH} is 
given as the following half-integral linear combination of cotangents of odd multiples of~$\pi/4n$:\
\be\label{hmnformula} h_{m,n} \= \sum_{k=1}^{2n-1} \ch(k)\,S\Bigl(\frac{km}n\Bigr)
  \,\cot\Bigl(\frac{\pi k}{4n}\Bigr) \;. \ee
\end{proposition}
\noindent{\it Proof.} Since the summand in~\eqref{HS} is an even function of~$k$, we can rewrite the sum
as a sum over all non-zero integral values of~$k$, and then use the periodicity (with period at most~$4n$)
of the factor $\ch(k)S(km/n)$ together with Euler's formula $\sum_{n\in\Z}(x+n)^{-1}=\pi\cot(\pi x)$ (with
the sum being interpreted as a Cauchy principal value) to get
\bas h_{m,n} &\= \frac{2n}\pi\,\sum_{k\in\Z\sm\{0\}} \frac{\ch(k)}k\,S\Bigl(\frac{km}n\Bigr) 
\= \frac{2n}\pi\,\sum_{\text{$k$ (mod~$4n$)},\;k\ne0} \ch(k)\,S\Bigl(\frac{km}n\Bigr)\,\sum_{r\in\Z}\frac1{4nr+k} 
\\  &\= \frac12\,\sum_{0<k<4n} \ch(k)\,S\Bigl(\frac{km}n\Bigr)\,\cot\Bigl(\frac{\pi k}{4n}\Bigr) 
\= \sum_{0<k<2n} \ch(k)\,S\Bigl(\frac{km}n\Bigr)\,\cot\Bigl(\frac{\pi k}{4n}\Bigr)  \eas
as  claimed.  \endpf
  
{\it Remark.} One can also deduce the formula for $c_{m,n}$  in Proposition~\ref{prop:selection}
from~\eqref{hmnformula}; this is left as an exercise to the reader.  We also observe 
that the terms on the right-hand side of~\eqref{hmnformula} (assuming, as we can, that 
$m$ and~$n$ are coprime) are integral multiples of cotangents of odd multiples of $\pi/4n$, and hence are
algebraic integers (in some cyclotomic field), except for the $k=n$ term if $m$ and $n$ are odd, which is 
half-integral because then $S(km/n)=\frac12$ and $\cot(\pi k/4n)=1$. Hence $h_{m,n}$ is an 
algebraic integer unless both $m$ and~$n$ are odd, in which case it equals an algebraic integer plus~$1/2$, and 
the sum $c_{m,n}=h_{m,n}+h_{n,m}$ is always an algebraic integer.

\medskip


\section{The connection with the generalized Riemann hypothesis} \label{sec:GRH}
In the first subsection of this section we outline the basic strategy for proving Theorem~\ref{GRH}.
The easier direction (a)$\,\Rightarrow\,$(b) of this theorem is then proved in Subsection~\ref{sec:PropABC} 
and the more difficult converse direction in Subsection~\ref{sec:Beurling}, using ideas from two classical papers 
by A.~Beurling. As already stated in the introduction, we will actually prove the following stronger result:
\begin{theorem} \label{Thm3} Let $L(s)$ be the $L$-series defined in~\eqref{chi4L4} and $R(N)$ be 
the function defined in~\eqref{defr}.
\begin{itemize} \item[(a)] If $L(s)$ has no zeros with $\Re(s)>\frac12$,
then $R(N)$ tends to infinity as $N$ tends to infinity.
\item[(b)] If $L(\r)=0$ for some $\r\in\CC$ with $\Re(\r)>\frac12$, then $R(N)$ satisfies the bound
 \be\label{bound} R(N) \;\le\; \frac\pi8\,\frac{|\r|^2}{\Re(\r)\m\h} \qquad(\forall N\ge1)\,. \ee
\end{itemize}  \end{theorem}

\subsection{The square-wave function~$\H$ and its multiplicative shifts} 
Set $\L^2 = \L^2(0, 1)$ and let $\xi$ be the characteristic function of the interval $(0, 1]$ and
$\H:\Rp\to\R$ the function
 \be \H(x) \= S(1/x) \= \begin{cases} 1 &\text{if $\frac1{4k+3}<x<\frac1{4k+1}$ for some $k\in\Z_{\ge0}$,}\\
       0 & \text{otherwise} \end{cases} \ee
(and $\H(x)=1/2$ if $1/x$ is an odd integer, but we are working in~$L^2$). Clearly both $\xi$ and~$\H$
 belong to~$\L^2$, since they are bounded with support in the interval $(0,1]$. On the other hand, if we
identify $L^2$ with the space of square-integrable functions on~$\Rp$ that vanish on~$(1,\infty)$, then
we have bounded operators $\t_a$ on $\L^2$ defined by  
$$ \t_af(x) \= f(x/a)  \qquad(0<a\le1,\;x\in\Rp) $$ 
We set $G_a=\t_aG$ and for any $\CD\subseteq\R$ and $N\in\N$ define subspaces $\V_\CD$ and $\V_N$ of~$L^2$ by
\be\label{defVs}  \V_\CD\,=\,\langle\, G_a\,|\,a\in\CD\cap(0,1]\,\rangle,
       \quad \V_N\,=\,\langle\, G_{n/N}\,|\,1\le n\le N\,\rangle\,, \ee
where $\langle\,S\,\rangle$ denotes the linear span of a subset~$S$ of a linear space. Finally, we denote by 
$\dst(\xi,\V_N)$ the $L^2$ distance from $\xi\in\L^2$ to the subspace $\V_N $, viz.
\be\label{defdist} \dst(\xi,\V_N) \= \bigl\|\,\xi\m\V_N\bigr\| 
       \= \inf_{x\in\R^N}\,\Bigl\|\,\xi\m \sum_{n=1}^Nx_n\,\t_{n/N}\H\,\Bigr\| \;. \ee
The key connection is then given by the following proposition.

\begin{proposition} \label{prA} The distance $\dst(\xi,\V_N)$ is related to the number $R(N)$ defined in~\eqref{defr} by
 \be \label{distance} \dst(\xi,\V_N)^2 \= \frac{\pi/4}{R(N)} \qquad(N\in\N)\,. \ee
 \end{proposition}
The proof of this, which is a straightforward computation, will be given in \S\ref{sec:PropABC} below.
On the other hand, a quite elementary argument, also given in \S\ref{sec:PropABC}, proves the following
\begin{proposition} \label{prB} If $L(\r)=0$ with $\Re(\r)>\h$ then
$\dst(\xi,\V_N)^2 \,\ge\, \frac{2\Re(\r)-1}{|\r|^2}$ for all~$N\in\N$ .
 \end{proposition}

\noindent Together these two results prove statement~(b) of Theorem~\ref{Thm3}. The harder statement~(a) 
follows from the two following results, proved in  \S\ref{sec:PropABC} and~\S\ref{sec:Beurling} respectively.

\begin{proposition} \label{prC} The space $\V_\CD$ is dense in $\V_\R$ for any set~$\CD$ 
that is dense in~$(0,1)$. \end{proposition}
 
\begin{proposition} \label{prD} If $L(s)$ has no zeros in~$\sg>\frac12$, then $\V_\R$  is dense in~$\L^2$. \end{proposition}

The whole argument is summarized by the following diagram, where the flow of the double arrows shows
that the six boxed statements are mutually equivalent.
\begin{center}
  \begin{tikzpicture}[scale=4,>=implies]
    \node (A) at (0,0) {$\boxed{\GR}$};
    \node (B) at (1.3,0) {$\boxed{\xi\in\overline{\V_\Q}}$};
    \node (C) at (2.8,0) {$\boxed{R(N) \text{ is unbounded}}$};
    \node (D) at (0,-0.86) {$\boxed{\overline{\V_\R}=L^2}$};
    \node (E) at (1.3,-0.86) {$\boxed{\overline{\V_\CD}=L^2\text{ ($\CD$ dense)}}$};
    \node (F) at (2.8,-0.86) {$\boxed{\lim\nolimits_{N\to\infty}R(N)=\infty}$};
    \draw[<-,double equal sign distance] (A) --node[above] {Prop.~\ref{prB}} (B);
    \draw[<-,double equal sign distance] (B) --node[above] {Prop.~\ref{prA}} (C);
    \draw[->,double equal sign distance] (A) --node[left]  {Prop.~\ref{prD}} (D);
    \draw[->,double equal sign distance] (E) --node[below] {Prop.~\ref{prA}} (F);
    \draw[->,double equal sign distance] (D) --node[below] {Prop.~\ref{prC}} (E);
    \draw[<-,double equal sign distance] (B) --node[right] {trivial} (E);
    \draw[<-,double equal sign distance] (C) --node[right] {trivial} (F);
  \end{tikzpicture}
\end{center}  
(The implications labelled ``Prop.~\ref{prA}" follow in one direction because each $\V_N$ is contained in~$\V_\Q$ 
and in the other because if $R(N)$ were bounded for some infinite sequence of integers $\{N_i\}$ then the set 
$\cup_i\V_{N_i}$ would be dense in~$\V_\R$ but have a non-zero distance from~$\xi$.)
\smallskip

\subsection{Proof of Propositions~\ref{prA}, \ref{prB} and~\ref{prC}} \label{sec:PropABC} 
Recall from linear algebra that the distance $d(v,V)$ from a vector~$v$ in a 
(real) Hilbert space to a finite-dimensional vector space~$V$ with basis $\{v_1,\dots,v_N\}$ is given by 
the formula $d(v,V)^2=\det(\wh G)/\det(G)$, where $G$ is the $N\times N$ Gram matrix with $(m,n)$-entry
$\la v_m,v_n\ra$ (scalar product of the basis elements $v_m$ and~$v_n$) and $\wh G$ the augmented
$(N+1)\times(N+1)$ matrix obtained from $G$ by adding a row and column with entries $\la v,v_n\ra$
and a diagonal entry~$\la v,v\ra$.  We apply this to the vector $v=\xi\in L^2$ and the space $\V_N$ 
defined as above as the span of the vectors $G_{n/N}\in L^2$ ($n=1,\dots,N$). (Here we mean the real 
rather than the complex span, i.e., we work in the real Hilbert space of real-valued square-integrable 
functions on~$(0,1)$.) We also denote by $\wV_N$ the span of $\V_N$ and~$\xi$.  We will show
in a moment that the vectors $G_{n/N}$ and $\xi$ are linearly independent, so 
$\dim\wV_N=\dim\V_N+1=N+1$. We can then compute the entries of the Gram matrices as
$$ \la\xi,\xi\ra \,=\,\int_0^1dx\,=\,1,\qquad \la\xi,G_{n/N}\ra
  \,=\,\int_0^1S\bigl(\frac n{Nx}\Bigr)\,dx\,=\,\frac nN\,L(1)\,=\,\frac{\pi n}{4N} $$
(the first is trivial and the second a special case of equation\eqref{scalprod} below) and
$$ \la G_{m/N},\,G_{n/N}\ra \,=\, \int_0^1 S\bigl(\frac m{Nx}\Bigr)\,S\bigl(\frac n{Nx}\Bigr)\,dx
 \= \frac mN\,C\Bigl(\frac mn\Bigr) \= \frac\pi{4N}\,c_{m,n}$$
(combine equations~\eqref{intC} and~\eqref{defC}). Proposition~\ref{prA} follows.

To find the connection with the $L$-series $L(s)$ and its complex zeros, we next
consider the ``complex monomial functions" $\mon_s$ on $(0, 1]$ defined  by 
$$\mon_s(t) \= t^{s-1} \xi(t) \= t^{s-1}    \qquad(0<t\le 1,\; s\in \mathbb{C})\,.$$ 
Note the key fact that $\mon_s \in {\L}^2$ if and only if $\sg > \frac12$, and also that for
 $f \in \L^2$ and $\sg > \frac12$, we have $\langle \mon_s, f \rangle $ = $\wt f(s)$,  
where $\wt{f}(s)$ denotes the Mellin transform  
$$  \wt f(s) \= \int_0^\infty f(x)\, x^{s-1}\,dx  \= \int_0^1 f(x)\, x^{s-1}\,dx\;. $$
(Here we are again identifying $\L^2=\L^2(0,1)$ with the space of functions $\L^2(\Rp)$ supported on~$(0,1]$.)
Then obviously 
\be\label{easy} \la\mon_s,\xi\ra \= \int_0^1 x^{s-1}\,dx \=\frac1s\,, \ee
and the Mellin transform of $G_a$ is almost equally easy to calculate using~\eqref{step}:
\be \label{scalprod} \la\mon_s,\H_a\ra\=\int_0^\infty S\Bigl(\frac ax\Bigr)\,x^{s-1}\,dx 
 \= \frac1s\,\sum_{k=0}^\infty\,\Bigl(\frac{a^s}{(4k+1)^s}\m\frac{a^s}{(4k+3)^s}\Bigr) 
  \= \frac{a^s}s\, L(s)\,, \ee
the calculation being valid for all $s$ in $\Re(s)>0$ since the integral then converges absolutely.
One consequence of this is that the functions $G_a$ ($0<a\le 1$) are all linearly independent, since
if some (finite) linear combination of them vanished then the vanishing of its Mellin transform would
imply that the product of $L(s)$ with a finite Dirichlet series vanished identically (first for~$\Re(s)>0$,
and then for all~$s$ by analytic continuation).  In particular, the space~$\V_N$ has dimension
exactly~$N$.  In fact, almost the same argument also shows that the functions $G_a$ and $\xi$ are
linearly independent, since a linear relation among them would imply that the product of $L(s)$ with 
a finite Dirichlet series was constant, which is obviously impossible.  It follows that the extended
space $\wV_N$ spanned by the functions $\H_{n/N}$ ($1\le n\le N$) and~$\xi$ has dimension~$N+1$,
as claimed above.

The proof of the direction (a)$\,\Rightarrow\,$(b) of Theorem~\ref{GRH} now follows easily.  
If there exists $\r$ with $ \re(\r) > \frac12$ and $\L(\r) = 0$  then by~\eqref{scalprod} we have  
$\mon_{\r} \perp \V_N $ for all $N$.  Since $\V_N$ is finite-dimensional, the infimum 
in equation~\eqref{defdist} is attained, so that $R(N)$ is unbounded (and hence $\dst(\xi,\V_N)$ has
lim inf equal to~0) if and only if there is a sequence of functions $h_N \in {\V}_N$  converging in 
$\L^2$ to~$\xi$; and this would contradict $\mon_{\r} \perp \V_N $ for all $N$ since $\la\mon_\r,\xi\ra=\r\iv\ne0$. 

This argument actually proves the stronger statement given in Proposition~\ref{prB} and in part~(b) of 
Theorem~\ref{Thm3}. Indeed, if $\r$ is a zero of $L(s)$ with $\Re(\r)>\h$, then on the one hand $\mon_\r$ belongs
to~$L^2$ as already noted, with $L^2$ norm
$$  \bigl\|\mon_\r\bigr\|^2 \= \int_0^1\bigl|x^{\r-1}\bigr|^2\,dx \= \frac1{2\,\Re(\r)\m1}\,, $$
and on the other hand $\la\mon_\r,f\ra=0$ for any $f\in\V_N$ by virtue of~\eqref{scalprod}. We apply this 
to the function $f=f_N\in\V_N$ attaining the minimum distance to~$\xi$.  Then applying the Cauchy-Schwarz 
inequality to the scalar product $\la\mon_\r,\xi-f_N\ra$ we get (using~\eqref{easy})
$$ \frac1{|\r|^2} \= \bigl|\la\mon_\r,\xi\ra\bigr|^2 \= \bigl|\la\mon_\r,\xi-f_N\ra\bigr|^2 
 \;\le\; \|\mon_\r\|^2\,\|\xi-f_N\|^2 \= \frac{d(\xi,\V_N)^2}{2\,\Re(\r)\m1}\,, $$
and~\eqref{bound} now follows immediately from~\eqref{distance}.

Finally, we prove the  Proposition~\ref{prC} (``density lemma").  We first note that the $L^2$ norm
$\|\H_a-\H_b\|$ tends to 0 as $a$ tends to~$b$, since by the calculation $\la\H_a,\H_b\ra=aC(a/b)$
given above (for rational numbers, but true in general by the same argument) we have
$$  \|\H_a-\H_b\|^2 \= \la\H_a,\H_a\ra\+\la\H_a,\H_b\ra\m2\,\la\H_a,\H_b\ra \= (a+b)\,C(1)\m 2a\,C(b/a)\,, $$
which tends to 0 as $b\to a$ because $C(x)$ is a continuous function. It follows immediately that if
$\CD$ is dense in~$(0,1)$ then the vector $\H_a\in L^2$ for any $a<1$ lies in the closure of $\V_\CD$,
and since $\V_\R$ is spanned by such vectors, this shows that $\overline{\V_\R}=\overline{\V_D}$ as claimed.

\subsection{Proof of Proposition~\ref{prD}} \label{sec:Beurling} 
It remains to show that $\GR$ implies that $\overline{\V_\R}=L^2$.

Let $\Hd^2 = \Hd^2(\CR) $ denote the Hardy space of the right half-plane $\CR=\{s\mid\Re(s)>\h\}$, 
that is, the space of all holomorphic functions $F(s)$ on $\CR$ that are square-integrable on the 
lines $\Re(s)=c$ for all $c>\h$ and with these $\L^2$ norms $\|\ \|_c$   bounded in $c$, and then 
equipped with the  Hilbert space norm $\|F\| = \sup \{\|F\|_c \mid c>\h\}$.  
We quote some results from the standard analysis of these spaces, citing as we go the relevant sections 
in the books by Hoffman~\cite{Ho} and Garnett~\cite{Gar} for the details. 

The Paley-Wiener theorem (\cite{Ho} pp.\,131--132) gives the fact that the Mellin transform is 
(up to scalar multiple) an isometry between the spaces $\L^2 $ and $\Hd^2$. Applying this to the 
function $\H$  implies that its Mellin transform $\,\wt\H(s)=\L(s)/s\,$  lies in the space $\Hd^2$, 
a fact which is also easy to see directly from the definition of $\Hd^2$. 
Next, we note that if we assume that $\GR$ holds, so that $\wt\H(s)$ has no zeros in~$\CR$, 
then the factorization theorem for functions in the space $\Hd^2$ (\cite{Ho} pp.\,132--133) applied
to the function $\wt\H(s)$  gives the implication 
\be\label{implication}
\text{\it $\L(s)$ has no zeros in $\sg>\tfrac12 \quad\Longrightarrow\quad \wt\H(s)$ is an outer function.}\ee
Indeed, the factorization theorem represents any element $F$ of $\Hd^2$ as the product of a ``Blaschke product,"
a ``singular function," and an ``outer function." Here we can omit the definitions of all three of these
since the Blaschke product is defined as a product over the zeros of $F$ in~$\CR$ and hence is equal to~1 
if $F$ has no such zeros, and the singular function factor is also constant because  $\wt\H(s)$ both continues 
across the line $\Re(s)=\h$ and has slow rate of convergence to zero as $s \to \infty$ on the positive reals 
(cf.~3.14 in~\cite{BS} and Chapter~II, Theorem~6.3, of~\cite{Gar}).  This leaves only the outer function 
factor of $\wt\H(s)$, and here we can apply a corollary of a theorem of Beurling on shift-invariant 
subspaces in $\Hd^2$  (\cite{Be1}, \cite{Ho} pp.\,99--101), which tells us that a function $F\in\Hd^2$ 
is outer if and only if the space 
$$ \CW(F)  \: \CC[1/s]\,F(s) \= \bigl\la\, s^{-n}F(s)\mid n\ge0\,\bigr\ra $$
is dense in $\Hd^2$. (Notice that $\CW(F)$ is contained in~$\Hd^2$ because the definition
of the $\Hd^2$ norm implies that $\|s\iv F(s)\|_{\Hd^2}\le 2\|F\|_{\Hd^2}<\infty$ for $F\in\Hd^2$,
so $s\iv\Hd^2\subset\Hd^2\,$.) 

It follows that GRH implies that $\CW(\wt\H)$ is dense in $\Hd^2$.
We can now apply the inverse Mellin transform from $\Hd^2$ to $\L^2$ to obtain the 
isometrically equivalent picture in~$\L^2$.  An easy calculation shows that $s\iv\wt f(s)$
is the Mellin transform of $If$ for any~$f\in L^2$, where 
$$ If(x) \: \int_x^\infty f(t) \, \frac{d t}{t} \;, $$
and hence that the space $\CW(\wt f)$ is the Mellin transform of the space
$$ W(f)  \: \CC[I]\,f \= \bigl\la\, I^n f\mid n\ge0\,\bigr\ra \qquad(f\in L^2).$$
It is not actually obvious that $I$ preserves the space~$L^2$, but this follows from the Paley-Wiener 
theorem together with the fact just shown that multiplication by $s\iv$ preserves~$\Hd^2$. A more direct 
argument is that the operator $I$ is the adjoint $\A^{\ast}$ of the averaging operator 
 $$\A h(x) \= \frac1x\,\int_0^x h(t)\,dt \,, $$
as is easily checked, and since $\A$ is a bounded operator on $\L^2(0,\infty)$ (cf.~\cite{Ha}, Example~5.4, p.\,23),
its adjoint~$I$ is also bounded on $\L^2(0,\infty)$.  But~$I$ also preserves the
condition of having support in~(0,1), so it is a bounded operator on~$L^2$.

As a consequence the implication~\eqref{implication} translates on the $L^2$ side to the statement
\be\text{\it $\L(s)$ has no zeros in $\sg>\tfrac12 \quad\Longrightarrow\quad W(\H)$ is dense in~$L^2$.}\ee
Proposition~\ref{prD} is then obtained by applying to $f=\H$ the following lemma, in which $\V(f)$ 
for any~$f\in L^2$ denotes the span of the multiplicative translates~$\t_af$~($0<a\le1$).

\smallskip
\noindent{\bf Lemma.} {\it For $f \in \L^2$ we have $W(f)\subset \overline{\V(f)}\,$.}

\smallskip\noindent{\it Proof.} We first note that it suffices just to show that $If \in \overline{\V(f)}$ for 
all~$f\in L^2$, because then $I^n f \in \overline{\V(f)}$ for all~$n\ge1$ by induction and the 
remark that $g\in\overline{\V(f)}\implies\overline{\V(g)}\subseteq\overline{\V(f)}$. Since 
$(W^\perp)^\perp = \overline{W}$ for all subspaces $W$, we need to show that 
\be \la h, \t_af \ra = 0 \;\, \text{for all $a\in(0,1]$}\, \implies \, \la If,\,h \ra = 0 .\ee
But since $a\t_{a\iv}$ is the adjoint of~$\t_a$ and $\A$ is the adjoint of~$I$, this statement just says
that any function that is orthogonal to all $h(ax)$ with $0<a\le 1$ is also orthogonal to $\A h$, and this is
obvious because $\A h(x)=\int_0^1h(ax)\,da\,$. \endpf

\medskip


\section{Analytic properties of $C(x)$ and $H(x)$} \label{sec:anal}

In this section we study the analytic properties of $C(x)$ and~$H(x)$, discussing in particular 
the question of their extendability from~$\Q$ to~$\R$ (which we answer completely only in 
the case of~$C(x)$) and describing the asymptotic properties in neighborhoods of rational 
points needed to establish Theorem~\ref{continuity} and its corollary (quantum modularity). 
We also discuss the relation of these two functions to the Dirichlet series~$L(s)$. It
should be mentioned in passing that much of the material in this section is similar in spirit to the
papers~\cite{BC1} and~\cite{BC2} of Bettin and Conrey, with $\ch$ replaced by the trivial character
and $L(s)$ by the Riemann zeta function.

\subsection{Sum and integral representations of $C(x)$ and $H(x)$}\label{Sec4.1}
Proposition~\ref{prop:analytic} below gives a number of formulas for these two functions.  
However, the word ``function" has a slighly different meaning in the two cases.  We have already 
defined them both as functions on~$\Q$. Part~(i) of the proposition says that $C$ extends 
continuously to~$\R$, and parts~(iv) and~(v) give further properties of that function.
But in the case of~$H$ we do not know to what extent it can be defined as a function
on~$\R$ (this question is discussed in Subsection~\ref{Sec4.2}, though not in great detail
since this is not our main interest).  However, it can certainly be defined as a distribution
on~$\R$, in fact even as the derivative of a continuous function, since by integrating~\eqref{HS} 
we get the absolutely convergent sum representation
$$\int_0^xH(x')\,dx'\=\frac{\pi x}8\m\frac12\sum_{k=1}^\infty\frac{\ch(k)}{k^2}\,A\Bigl(\frac{kx}2\Bigr) $$
for its integral, where $A$ is the periodic continuous function on~$\R$ defined by $A(n+\ep)=(-1)^n\ep$
for $n\in\Z$ and $|\ep|\le\frac12$. Parts (ii) and~(iii) of the proposition are then to be interepreted in the
sense of distributions, one giving the Fourier expansion of the even periodic distribution~$H$ and
the other its relation, already mentioned in the introduction, with the derivative of~$C$.

\begin{proposition} \label{prop:analytic}
$($i$)$ The function $C(x)$ extends continuously from~$\Q$ to~$\R$ via the formula 
\be\label{intC}  C(x) \= \int_0^\infty S(t)\,S(xt)\,\frac{dt}{t^2}\qquad(x\in\R)\,.\ee

$($ii$)$ The distribution $H$ has the Fourier expansion
\be\label{HFour} H(x)\=\frac\pi8\m\frac2\pi\,\sum_{n=1}^\infty
       \frac{\ch(n)\,d(n)}n\,\cos\Bigl(\frac{\pi nx}2\Bigr)\,,\ee
where $d(n)$ denotes the number of positive divisors of~$n$.

$($iii$)$ The function $C(x)$ and the distribution $H(x)$ are related by the equation
  \be\label{Cderiv}  C'(x) \= H(1/x)\qquad(x>0)\,. \ee

$($iv$)$ The function $C(x)$ is also given by the sum
\be\label{CbyJ} C(x)\=\frac\pi8 \+ x\,\sum_{n=1}^\infty\ch(n)\,d(n)\,J\Bigl(\frac{\pi nx}2\Bigr)\qquad(x>0)\,,\ee
where $J(x)$ denotes the modified cosine integral
\be\label{Jdef} J(x) \= -\int_x^\infty\frac{\cos(t)}{t^2}\,dt\qquad(x\in\Rp)\,. \ee

$($v$)$ The Mellin transform of $C(x)$ is given by
\be\label{Ctilde} \wt C(s) \= -\,\frac{L(-s)\,L(s+1)}{s\,(s+1)}  \qquad(-1<\s<0),\ee
where $L(s)$ is the analytic continuation of the $L$-series defined in~\eqref{chi4L4}.

\end{proposition}  

\noindent{\it Proof.} Formula~\eqref{intC} for $x\in\Qp$ is precisely the statement of
Proposition~\ref{prop:cmnintegral} after an obvious change of variables, and the convergence and continuity 
of the integral for all~$x\in\Rp$ are easy consequences of the facts that $S$ is supported on $[1,\infty)$ 
and is locally constant and bounded. To prove~(ii), we insert the (standard) Fourier expansion
\be\label{SFour} S(x)\=\frac12 \m 
  \frac2\pi \sum_{m=1}^\infty \frac{\ch(m)}m\,\cos\Bigl(\frac{\pi mx}2\Bigr)\,, \ee
 of the even periodic function $S(x)$ into the sum~\eqref{HS}, combine the double sum into a
single one, and use Leibniz's formula $L(1)=\pi/4$. For (iii), we
insert~\eqref{defS} into~\eqref{intC} to get  
 \be\label{Csum} C(x) \= \int_1^\infty \sum_{0<k<xt}\ch(k)\,S(t)\,\frac{dt}{t^2} 
\= \sum_{k=1}^\infty \ch(k)\,\int_{k/x}^\infty \frac{S(t)}{t^2} \,dt \,, \ee
from which equation~\eqref{Cderiv} follows immediately by differentiating and using~\eqref{HS}. (The easy 
justification of these steps in the sense of distributions by integrating against smooth test functions 
is left to the reader.) Finally, to prove~(iv) we combine~\eqref{Cderiv} and~\eqref{defhmn} to get
$$ -\,x^2\,\frac d{dx}\biggl(\frac{C(x)}x\biggr) \= C(x) \m x\,C'(x) 
  \= C(x) \m x\,H\Bigl(\frac1x\Bigr) \= H(x) $$
and then insert the Fourier development~\eqref{HFour} and integrate term-by-term. This 
gives~\eqref{CbyJ} up to an additive term $\lambda x$ that can be eliminated by noting that
$C$~is bounded. Note that the sum in~\eqref{CbyJ} is uniformly and absolutely convergent since
a simple integration-by-parts argument shows that $J(x)=\text O(x^{-2})$ as~$x\to\infty$. Finally, for~(v) 
we note first that the Mellin transform $\wt S(s):=\int_0^\infty S(x)\,x^{s-1}dx$ of~$S$ is given 
according to~equation~\eqref{scalprod} by $\wt S(s)=-L(-s)/s$ for $\s=\Re(s)<0$, and then use~\eqref{intC} to get
$$ \wt C(s) \= \int_0^\infty S(t)\,\frac{\wt S(s)}{t^s}\,\frac{dt}{t^2} \= \wt S(s)\,\wt S(-s-1)
 \= -\,\frac{L(-s)\,L(s+1)}{s\,(s+1)} $$
for $-1<\s<0$, as claimed. (The Mellin transform exists in this strip because $C(x)=\text O(\min(|x|,1))\,$
as an easy consequence of~\eqref{intC}.) In particular, $\wt C(s)$ extends meromorphically to all of~$\CC$ 
and is invariant under $s\mapsto-s-1$. We also note that by the well-known functional equation of~$L(s)$, 
equation~\eqref{Ctilde} can be written in the alternative form 
$$ \wt C(s)\= -\,\frac{(2/\pi)^{s+1}\G(s)\,\cos(\pi s/2)}{s+1}\,L(s+1)^2\,.  $$
We can then use this to give a second derivation of~\eqref{CbyJ}, not making use of the distribution~$H$,
by first applying the Mellin inversion formula to write $C(x)$ as an integral over a vertical line 
$\s=c$ with $-1<c<0$ and then shifting the path of integration to the right, picking up the term $\pi/8$
in~\eqref{CbyJ} from the residue at $s=0$ and permitting us to use the convergent Dirichlet series 
representation $L(s)^2 = \sum_{n=1}^\infty\ch(n)d(n)n^{-s}$ for~$\s>1$.  (We omit the details
of this calculation.)  \endpf

\subsection{The function $H$ on the real line}\label{Sec4.2}
In this subsection we give two one-parameter generalizations of $H(x)$ and discuss possible ways
to define this function at irrational arguments.

We first generalize the definitions of $H(x)$ and $C(x)$ and the assertions of Proposition~\ref{prop:analytic}
to families depending on a complex parameter~$s$. Specifically, we can generalize~\eqref{HS} to 
\be\label{defHs} H_s(x) \: \sum_{k=1}^\infty \frac{\ch(k)}{k^s}\,S(kx)\qquad(x\in\R,\;\s=\Re(s)>1)\,,   \ee 
which now converges absolutely. The same calculation as for $H(x)$ then gives the Fourier expansion
$$ H_s(x)\=\frac{L(s)}2\m\frac2\pi\,\sum_{n=1}^\infty
   \frac{\ch(n)\,\s_{s-1}(n)}{n^s}\,\cos\Bigl(\frac{\pi nx}2\Bigr)\,,$$
where $L(s)$ is the Dirichlet $L$-function defined in~\eqref{chi4L4} and
$\s_\nu$ is the sum-of-divisors function $\s_\nu(n)=\sum_{k|n,k>0}k^\nu$. Similarly,
we define $C_s(x)=x^sC_s(1/x)$ for all $s\in\CC$ with $\s>0$ by
\be\label{defCs} C_s(x) \= s\int_0^\infty\!S(t)\,S(xt)\,\frac{dt}{t^{s+1}}
  \= \sum_{j,\,k>0} \frac{\ch(j)\,\ch(k)}{\max(j/x,k)^s}\,,  \ee
where the equality of the two expressions is proved by the same calculation as for the 
proof of~\eqref{intC}.  Then $C_s(x)$ and $H_s(x)$ are related by 
\be\label{CsHs}   C_s(x) \= H_s(x)\+x^s\,H_s\Bigl(\frac1x\Bigr)\,, 
      \qquad C_s'(x) \= s\,x^{s-1}H_s\Bigl(\frac1x\Bigr)\qquad(x>0) \ee
by the same calculations as in the special case~$s=1$ and one has an expansion like~\eqref{CbyJ} but with 
$J(x)$ replaced by $-\int_x^\infty t^{-s-1}\cos(t)\,dt$.  As long as $\s>1$ all of the sums and integrals
involved are absolutely convergent and all steps are justified.

A different possible way to regularize $H(x)$ for $x\in\R$ is to define
\be\label{TepsDef} H(x) = \lim\limits_{\v\searrow0}\Bigl(\dfrac8\pi-\dfrac2\pi T(x,\v)\Bigr), \;\,\text{where}\;\,
 T(x,\v) \= \sum_{n=1}^\infty \frac{\chi(n)\,d(n)}n\,e^{-n\v}\,\cos(n\pi x/2)\,. \ee
Using the identity
 $\sum\limits_{n>0}\dfrac{d(n)}n\,x^n = \sum\limits_{k>0}\dfrac1k\,\log\Bigl(\dfrac1{1-x^k}\Bigr)$ 
and standard trigonometric identities, we find after a short calculation that
 $$ T(x,\v) \= \frac12\,\sum_{k=1}^\infty\frac{\chi(k)}k\,\arctan\biggl(\frac{\cos(k\pi x/2)}{\sinh(k\v)}\biggr)\,.$$
Both this series and the original one defining~$T$ converge exponentially fast for any positive~$\v$,
giving us another possible approach to the analytic properties of the limiting function~$H$.

Summarizing this discussion, we have at least five potential definitions of $H(x)$ for $x\in\R$:

\smallskip\noindent{\bf Definition 1}: Define $H(x)$ by the series~\eqref{HS}, if this sum converges.

\smallskip\noindent{\bf Definition 2}: Define $H(x)$ by the Fourier series~\eqref{HFour}, if this series converges.

\smallskip\noindent{\bf Definition 3}: Define $H(x)$ as the limit for $s\searrow1$ of the series~\eqref{defHs}, 
if this limit exists.

\smallskip\noindent{\bf Definition 4}: Define $H(x)$ as the limit in~\eqref{TepsDef}, if this limit exists.

\smallskip\noindent{\bf Definition 5}: Define $H(x)$ as the limit of $H(x')$ as $x'\in\Q$ tends to~$x$, 
if this limit exists.

\smallskip\noindent We can then ask---but have not been able to answer---the question whether any or all 
of these definitions converge for irrational values of~$x$ or, if they do, whether they give the same 
value.  The last definition is in the sense the strongest one, since any definition of $H$ 
on~$\R$ or a subset of~$\R$ that agrees with the original definition on~$\Q$ must coincide
with the value in Definition~5 at any argument~$x$ at which this function is continuous.
In any case, we pose the explicit question:

\smallskip\noindent{\bf Question}: Are the five definitions given above convergent and equal to one another
for all irrational values of~$x$, or for all $x$ belonging to some explicit set of measure~1\,?

\smallskip\noindent It seems reasonable to expect the answer to the second question to be affirmative
with the set of measure~1 being the complement of some set of irrational numbers having extremely good
rational approximations, like the well-known ``Brjuno numbers" in dynamics.

\subsection{Asymptotic behavior of $H(x)$ and $C(x)$ near rational points} \label{Sec4.3}
The results of Section~\ref{Sec4.1} prove the continuity of~$C(x)$, which is part of the 
statement of Theorem~\ref{continuity} in the Introduction.  We now discuss the remaining 
statements there, concerning the behavior of~$C$ and~$H$ near rational points. 

We begin by looking numerically at the asymptotic properties of $C(x)$ near~$x=1$.
Computing the values of $C\bigl(1\pm\frac1n\bigr)$ for $1\le n\le1000$ and using a numerical
interpolation technique that is explained elsewhere (see e.g.~\cite{GM}), we find
empirically an expansion of the form
\be \label{Cat1} C\Bigl(1\pm\frac1n\Bigr) \;\sim\; c_0 \+ \frac1n\,\Bigl(-\frac14\,\log n+c_1^\pm\Bigr) 
 \,\pm\,\frac{c_2}{n^2}\+\frac{c_3}{n^3} \,\pm\,\frac{c_4}{n^4}\+\frac{c_5}{n^5} \pm \cdots \ee
with the first few coefficients having numerical values given by
\bas & c_0=\,C(1)=\frac\pi4\,,\quad c_1^+=c_1^-+c_0=-0.23528274848426799887\cdots\,, \quad c_2\,=\,-\frac18\,,
 \\   & c_3\,=\,0.058801396529669\cdots\,,\quad c_4\,=\,\frac1{48}\m c_3\,,\quad c_5 = 0.01927655829\cdots\;. \eas
These numerics become clearer if we work with $H(x)$ instead, where we find the simpler expansion
\be \label{Hat1} H\Bigl(1\pm\frac1n\Bigr) \;\sim\; 
 \pm\,\biggl(\frac14\,\log n \+ h_0^\pm \+ \frac{h_2}{n^2}\+\frac{h_4}{n^4} \+ \cdots \biggr) \ee
with no odd powers of $1/n$ and with coefficients given numerically by
$$ h_0^+=h_0^-+c_0=0.7706809118817\cdots\,, \quad h_2=0.01713472986\cdots\,,
       \quad h_4=-0.00345272385\cdots\,, $$
which then gives the expansion~\eqref{Cat1} with $c_1^\pm=-h_0^\mp-\dfrac14$, 
$c_3=h_2+\dfrac1{24}$, and more generally $c_r=(-1)r\sum_{0<i<r/2}\binom{r-1}{2i-1}h_{2i}$
for $r>1$. Moreover, on calculating the next few
values of $h_{2i}$ numerically to high precision we are able to guess the closed formulas
\be\label{Bern}  h_0^\pm \= \pm\,\frac\pi8\+\frac\gamma4 \+ \frac14\log\frac8\pi\,, \qquad
h_{2i} \= (-1)^{i-1}\,\frac{(2^{2i-1}-1)^2\,B_{2i}^2}{2i\,(2i)!}\,\pi^{2i}\;\quad(i\ge1) \ee
for the coefficients in~\eqref{Hat1}, where $\g$ is Euler's constant and $B_n$ the $n$th 
Bernoulli number. In fact, it is not difficult to prove~\eqref{Hat1} (and hence also~\eqref{Cat1})
with these values of~$h_i$ by applying a twisted Euler-Maclaurin formula (giving the asymptotics 
as $\v\to0$ of sums over an interval of $\ch(n)f(n\v)$ for a smooth function~$f$) to the closed formula
$$   H\Bigl(1\pm\frac1n\Bigr) \= \frac{4n}\pi\,\sum_{0<k<2n\atop k\equiv1\!\!\pmod4}\cot\frac{\pi k}{4n} $$
which follows easily from~\eqref{hmnformula}, though we will not carry this out here.

The ``$\log n$" terms in equations~\eqref{Cat1} and ~\eqref{Hat1} already show that $C(x)$ is not
differentiable at $x=1$ and $H(x)$ is not continuous at~$x=1$, but from these formulas one might imagine 
that the functions $C(x)-|x-1|\log(|x-1|)/4$ and $H(x)-\log(|x-1|)/4$ are $C^\infty$ to both the right
and the left of this point.  However, this is not the case, because both equations~\eqref{Cat1} and~\eqref{Hat1} 
are valid only when $n$ is an integer and change in other cases.  For instance, if $n$ tends to infinity in
$\N+\frac12$ rather than~$\N$, then $H(1\pm1/n)$ has an expansion of the same form as~\eqref{Hat1} and 
with the same constants $h_0^\pm$, but with $h_2=\pi^2/576$ replaced by $-7\pi^2/1152$, $h_4=-49\pi^4/1382400$
replaced by $127\pi^4/2764800$, etc.) This statement, which again can be proved using an appropriate 
twisted version of the Euler-Maclaurin formula, is a typical phenomenon of quantum modular forms.  
 
If we look at the asymptotics of $C(x)$ and $H(x)$ as $x$ tends to any rational number~$\alpha$ with odd
numerator and denominator, then we find a similar behavior, which is not surprising since any such~$\alpha$ 
is $\G_\th$-equivalent to~1 and both $H$ and $C$ have transformation properties, modulo functions with 
better smoothness properties, under the action of~$\G_\th$.  More precisely, if we write $\a=a/c$ with 
$a$ and $c>0$ odd and coprime and complete $\bigl(\smallmatrix a\\c\endsmallmatrix\bigr)$ 
to a matrix $\sma abcd\in\SL$, then we find asymptotic expansions
 \bas H\bigl(\frac{an+b}{cn+d}\bigr) &\;\sim\; 
      \pm\,\frac{\log|cn+d|}{4c} \+ \sum_{i=0}^\infty\frac{h^\pm_{i,\a}}{(cn+d)^i}\,, \\
  C\bigl(\frac{an+b}{cn+d}\bigr) &\;\sim\; C\bigl(\frac ac\bigr)
    \,\pm\,\frac{\log|cn+d|}{4(cn+d)} \+ \sum_{i=1}^\infty\frac{c^\pm_{i,\a}}{(cn+d)^i} \eas
as $n\to\pm\infty$ with $n\in\Z$, and similar expansions with other coefficients $h^\pm_{i,\a,\b}$
and $c^\pm_{i,\a,\b}$ if $n\to\pm\infty$ with $n\in\Z+\b$ for some fixed rational number~$\b$.

Similar statements hold for $x$ tending to rational numbers~$\a$ with an even numerator or denominator,
but now without the log term, the simplest case being
 $$ C\bigl(\frac1n\bigr) \;\sim\; \frac{\pi}{8n} \+ \ch(n)\,\frac{\pi^2}{16n^2} \+ 
  \ch(n-1)\,\frac{\pi^3}{64n^3} \+ \cdots 
 \+ \frac{\ch(n-k+1)}4\,\frac{A_k^2}{k!}\,\Bigl(\frac\pi{2n}\Bigr)^{k+1} + \cdots $$
where $A_k$ is the number of ``up-down" permutations of~$\{1,\dots,k\}$ (permutations~$\pi$ where
$\pi(i)-\pi(i-1)$ has sign $(-1)^i$ for all~$i$), which is also equal to the coefficient
of $x^k/k!$ in $\tan x+\sec x$. (Note that $A_{k-1}=(4^k-2^k)|B_k|/k$ for $k$~even,
so that the coefficients here are closely related to those in~\eqref{Bern}.)  The difference
between the two cases is due to the fact that the action of $\G_\th$ on $\mathbb P^1(\Q)$ has two orbits 
(cusps), and comes from the relation that we will see in the final section between the quantum modularity
properties of~$C$ and $H$ and a specific Maass modular form on~$\G_\th$ that is exponentially small
at the cusp~$\infty$ but has logarithmic growth at the cusp~1\,.


\section{Generalization to other odd Dirichlet characters} \label{sec:OtherChi}

In this section we discuss the changes that are needed when the character defined
in~\eqref{chi4L4} is replaced by an arbitrary odd primitive Dirichlet character. Suppose that
$D<0$ is the discriminant of an imaginary quadratic field~$K$ (or equivalently, that $D$ is
either square-free and congruent to 1~mod~4 or else equal to 4 times a square-free number not
congruent to~1 mod~4), and let $\ch=\ch_D=\bigl(\frac D\cdot\bigr)$ be the associated Dirichlet 
character, the case studied up to now corresponding to $K=\Q(i)$ and~$D=-4$.  The Dirichlet 
$L$-series $L(s,\ch)$ is defined as $\sum_{n>0}\ch(n)n^{-s}$ for~$\sg>1$ and by analytic continuation 
otherwise, and is equal to the quotient of the Dedekind zeta function of~$K$ by the Riemann
zeta function. Its value at $s=1$ is well-known to be non-zero and given by
$$ L(1,\ch) \= \frac\pi{\sqrt{|D|}}\,h'(D)\,, $$
where $h'(D)$ is $1/3$ or $1/2$ if $D=-3$ or $D=-4$ and otherwise is the class number of~$K$. 

We now define $S(x)$ by the same formula~\eqref{defS} as before, with $\ch=\ch_D$.  It is still 
periodic (with period~$|D|$) and even and hence bounded, but now has average value $h'(D)$ 
rather than~1/2 and also no longer takes on only the values 0~and~1, as one sees in the
two following pictures of the graphs of this function for $D=-3$ and~$D=-7$.  
\medskip
\begin{center}
  \begin{tikzpicture}[scale=1.8]
    \newcommand{\charthree}{
      \path[fill=black!10!white] (1,0) rectangle (2,1);
      \draw[very thick] (0,0) -- (1,0);
      \draw[dashed] (1,0) -- (1,1);
      \draw[very thick] (1,1) -- (2,1);
      \draw[dashed] (2,1) -- (2,0);
      \draw[very thick] (2,0) -- (3,0);
      \path[draw,fill=white] (1,0) circle [radius=1pt];
      \path[draw,fill=white] (2,0) circle [radius=1pt];
      \path[draw,fill=white] (1,1) circle [radius=1pt];
      \path[draw,fill=white] (2,1) circle [radius=1pt];
      \path[draw,fill=black] (1,0.5) circle [radius=1pt];
      \path[draw,fill=black] (2,0.5) circle [radius=1pt];
    }
    \begin{scope}
      \clip(-1.5,-0.5) rectangle (5.5,1.5);
      \begin{scope}[xshift=-3cm]
        \charthree
      \end{scope}
      \begin{scope}[xshift=0cm]
        \charthree
      \end{scope}
      \begin{scope}[xshift=3cm]
        \charthree
      \end{scope}
    \end{scope}
    \draw[->] (-1.5,0) -- (5.6,0);
    \draw[->] (0,-0.1) -- (0,1.4);
    \foreach \x in {-1,0,1,...,5} {
      \draw (\x,0.05) -- (\x,-0.05) node[anchor=north] {$\x$};
    }
    \draw (0.05,0.5) -- (-0.05,0.5); 
    \draw (0.05,1) -- (-0.05,1) node[anchor=east] {$1$};
  \end{tikzpicture}\\[2mm]
\smallskip
  \begin{tikzpicture}[scale=1.3]
    \newcommand{\charseven}{
      \path[fill=black!10!white] (0,0) 
      \foreach \dy in {0,1,1,-1,1,-1,-1} { -- ++(0,\dy)  -- ++(1,0)} -- cycle;
      \draw[dashed] (0,0) 
      \foreach \dy in {0,1,1,-1,1,-1,-1} { -- ++(0,\dy)  -- ++(1,0)};
      \draw[very thick] (0,0) -- (1,0);
      \draw[very thick] (1,1) -- (2,1);
      \draw[very thick] (2,2) -- (3,2);
      \draw[very thick] (3,1) -- (4,1);
      \draw[very thick] (4,2) -- (5,2);
      \draw[very thick] (5,1) -- (6,1);
      \draw[very thick] (6,0) -- (7,0);
      \path[draw,fill=white] (1,0) circle [radius=1pt];
      \path[draw,fill=white] (1,1) circle [radius=1pt];
      \path[draw,fill=white] (2,1) circle [radius=1pt];
      \path[draw,fill=white] (2,2) circle [radius=1pt];
      \path[draw,fill=white] (3,1) circle [radius=1pt];
      \path[draw,fill=white] (3,2) circle [radius=1pt];
      \path[draw,fill=white] (4,1) circle [radius=1pt];
      \path[draw,fill=white] (4,2) circle [radius=1pt];
      \path[draw,fill=white] (5,1) circle [radius=1pt];
      \path[draw,fill=white] (5,2) circle [radius=1pt];
      \path[draw,fill=white] (6,0) circle [radius=1pt];
      \path[draw,fill=white] (6,1) circle [radius=1pt];
      \path[draw,fill=black] (1,0.5) circle [radius=1pt];
      \path[draw,fill=black] (2,1.5) circle [radius=1pt];
      \path[draw,fill=black] (3,1.5) circle [radius=1pt];
      \path[draw,fill=black] (4,1.5) circle [radius=1pt];
      \path[draw,fill=black] (5,1.5) circle [radius=1pt];
      \path[draw,fill=black] (6,0.5) circle [radius=1pt];
    }
    \begin{scope}
      \clip(-1.5,-0.5) rectangle (8.5,2.5);
      \begin{scope}[xshift=-7cm]
        \charseven
      \end{scope}
      \begin{scope}[xshift=0cm]
        \charseven
      \end{scope}
      \begin{scope}[xshift=7cm]
        \charseven
      \end{scope}
    \end{scope}
    \draw[->] (-1.5,0) -- (8.6,0);
    \draw[->] (0,-0.1) -- (0,2.4);
    \foreach \x in {-1,0,1,...,8} {
      \draw (\x,0.05) -- (\x,-0.05) node[anchor=north] {$\x$};
    }
    \draw (0.05,0.5) -- (-0.05,0.5); 
    \draw (0.05,1) -- (-0.05,1) node[anchor=east] {$1$};
    \draw (0.05,1.5) -- (-0.05,1.5); 
    \draw (0.05,2) -- (-0.05,2) node[anchor=east] {$2$};
  \end{tikzpicture}\\[2mm]
  \textbf{Figure 5. Graphs of the functions $S(x)$ for $D=-3$ and $D=-7$}
\end{center}

\medskip\noindent 
We define $H(x)=H_D(x)$ and $C(x)=C_D(x)$ for~$x\in\Q$ by the same formulas~\eqref{HS} and~\eqref{CH} as before,
but we define $h_{m,n}$ and $c_{m,n}$ by~\eqref{defH} and~\eqref{defC} with a factor $|D|/\pi$ instead of $4/\pi$ 
in order to get integral linear combinations of cotangents. 
The Fourier expansions~\eqref{SFour} and~\eqref{HFour} are then replaced by 
\begin{align}\label{newSFour} S(x) &\= h'(D) \m \frac{\sqrt{|D|}}\pi \sum_{m=1}^\infty 
   \frac{\ch(m)}m\,\cos\Bigl(\frac{2\pi mx}{|D|}\Bigr)\,, \\
  \label{newHFour}  H(x) &\= \frac{\pi\,h'(D)^2}{\sqrt{|D|}} \m \frac{\sqrt{|D|}}\pi\,
    \sum_{n=1}^\infty \frac{\ch(n)\,d(n)} n\,\cos\Bigl(\frac{2\pi nx}{|D|}\Bigr)\,, \end{align}
and the formula corresponding to~\eqref{hmnformula} now reads
\be\label{newhmn} h_{m,n}\=\sum_{0<k<|D|n/2} \ch(k)\,S\Bigl(\frac{km}n\Bigr)\,\cot\Bigl(\frac{\pi k}{|D|n}\Bigr)\ee
by the same calculations as before. Finally we also define the matrices $\C_N$ and $\hC_N$ exactly as we did
in the special case~$D=-4$ (except for replacing the lower right-hand entry in $\hC_N$ in~\eqref{defCN} by
$|D|N/\pi$), and then from the scalar product calculations
$$ \la\xi,\xi\ra\,=\,1\,,\qquad \la\xi,\H_a\ra\,=\,L(1)\,a\,,\qquad \la\H_a,\H_b\ra\,=\,a\,C(a/b)\,,$$
which are proved exactly as before, we deduce the same connection as in Theorems~\ref{GRH}
and~\ref{Thm3} between the unboundedness of the function~$R(N)$ and the Riemann hypothesis for~$L(s,\chi)$.

The only real difference with the case $D=-4$ is in the argument for quantum modularity.  The function 
$H(x)$ on~$\Q$ is still even and periodic of period~$|D|$ (and also anti-periodic up to a constant with 
period~$|D|/2$ if~$D$ is even), it again has discontinuities at infinitely many rational points by an
argument similar to the one given in~\ref{Sec4.3}, and by the analog of Proposition~\ref{prop:analytic} the 
function $C:\Q\to\R$ again has a continuous extension to~$\R$ and is therefore much better behaved analytically 
than $H(x)$, as illustrated by the following graphs of these functions for~$D=-3$, which look qualitatively
much like their $D=-4$ counterparts in Figures~2 and~3. 
\medskip
\begin{center}
  \includegraphics[width=0.8\linewidth]{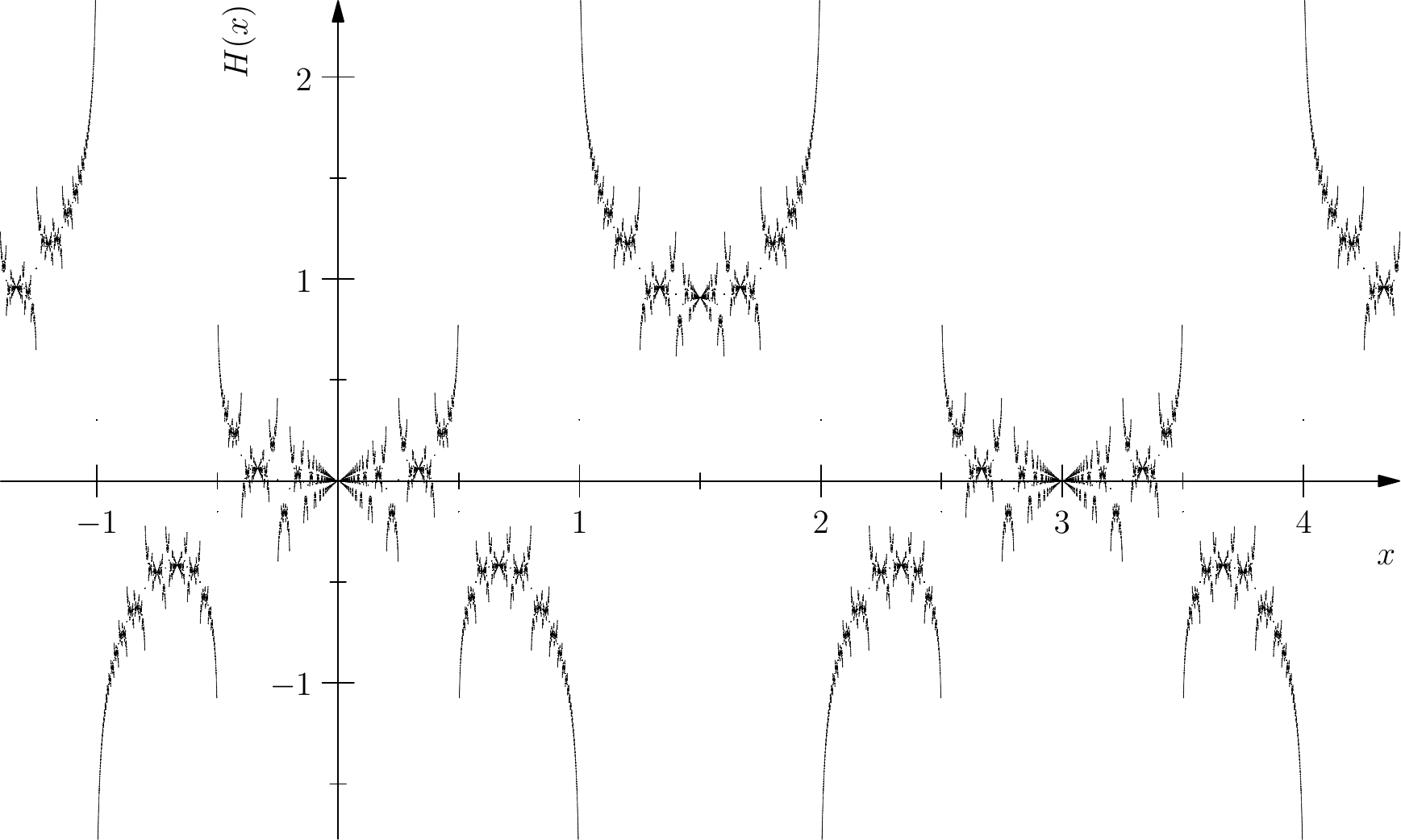}\\*[2mm]
  \textbf{Figure 6. Graph of the function $H(x)$ for $D=-3$}
\end{center}

\smallskip
\begin{center}
  \includegraphics[width=0.8\linewidth]{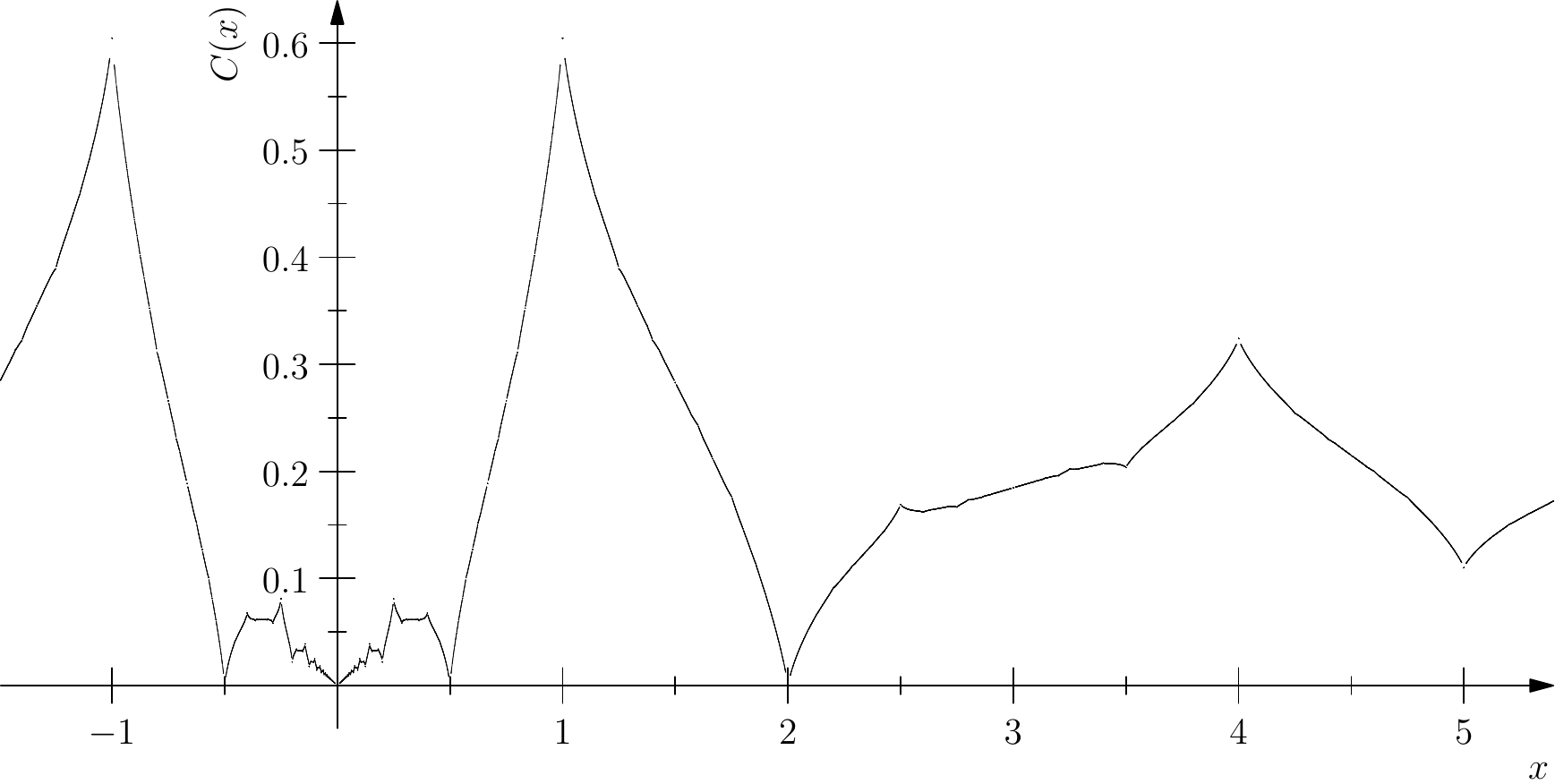}\\*[2mm]
  \textbf{Figure 7. Graph of the function $C(x)$ for $D=-3$}
\end{center}

\noindent The difference is that this longer suffices to prove the quantum modularity of~$H$ because the subgroup of~$\SL$ 
generated by the matrices $T^{|D|}$ (or $T^{|D|/2}$ if $D$ is even) and~$S$ now no longer has finite index. 
Instead, we need a statement like the corollary to Theorem~\ref{continuity} with $\G_\th$ replaced 
by the congruence subgroup $\G_{(D)}=\G_0(D)\cup S\G_0(D)$ of~$\SL$ (or by $\G_{(D/2)}$ if~$D$ 
is even, which is indeed $\G_\th$ if~$D=-4$).  This statement is given in the following theorem. 

\begin{theorem} \label{QMFforD} Let $D<0$ and $H(x)=\sum_{k>0}\ch_D(k)S_D(kx)/k$ be as above. Then the function
 \be\label{defCgam} C_\g(x) \: H(x) \m \v(\g)\,|cx+d|\,H\Bigl(\frac{ax+b}{cx+d}\Bigr) \qquad(x\in\Q) \ee
extends continuously to~$\R$ for all matrices $\g=\sma abcd\in\G_{(D)}$, where $\v:\G_{(D)}\to\{\pm1\}$
is the homomorphism mapping $\G(D)$ to~$1$ and $S$~to~$-1$.
\end{theorem}  

\noindent This theorem is illustrated for $D=-3$ and a typical element of~$\G_{(-3)}$ in the following figure. \break
\smallskip
\begin{center}  
  \includegraphics[width=0.8\linewidth]{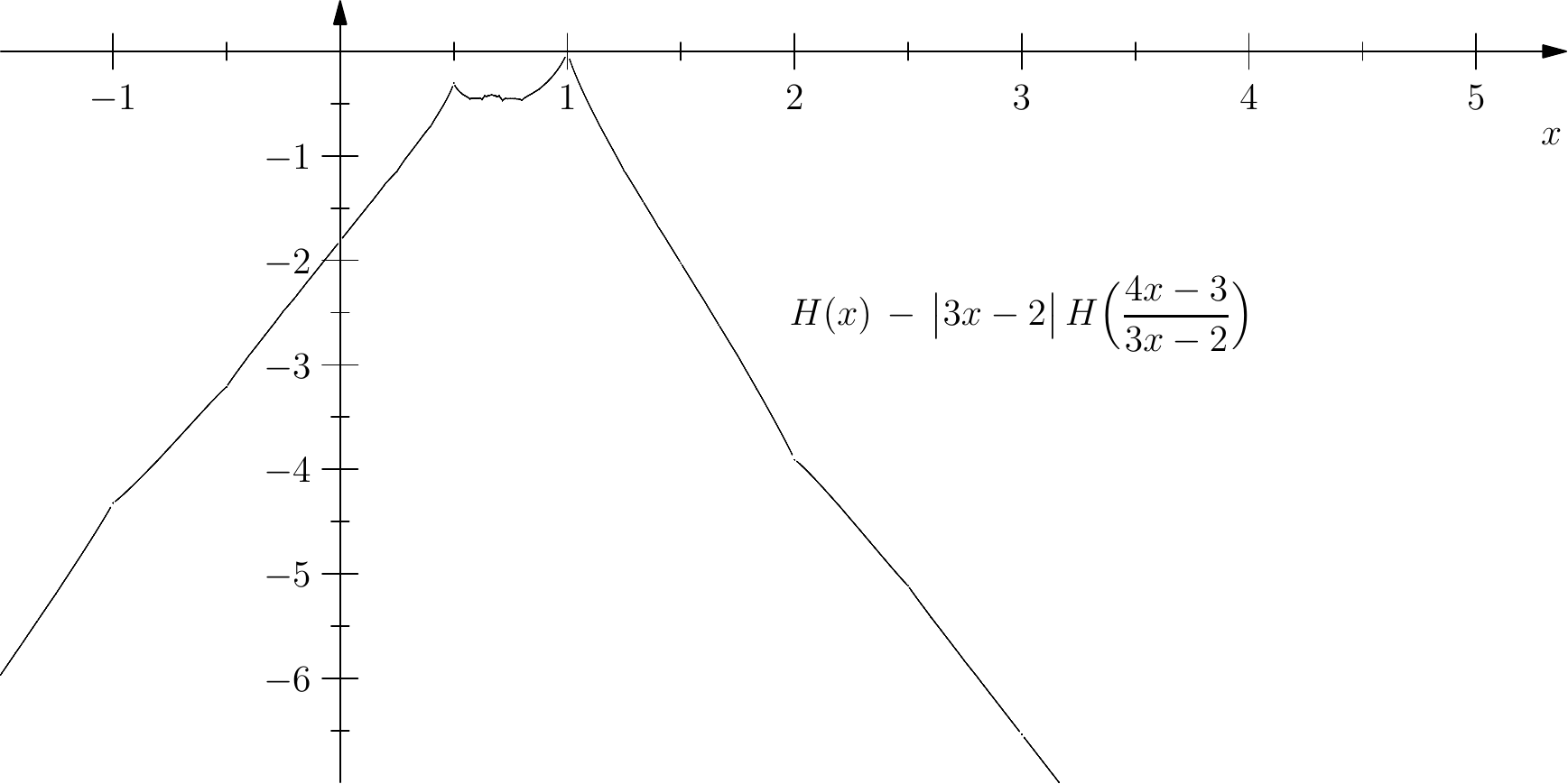}\\*[2mm]
  \textbf{Figure 8. Graph of $C_{\gamma}(x)$ for $D=-3$, $\gamma=\sma4{-3}3{-2}$}
\end{center}

\medskip

We will deduce Theorem~\ref{QMFforD} as a consequence of the following proposition,
which is a generalization to arbitrary elements $\g\in\G_{(D)}$ of equation~\eqref{Csum} for $D=-4$ and~$\g=S$.

\begin{proposition} \label{prop:Cgamint}
Let $\g=\begin{pmatrix} a&b\\c&d\end{pmatrix}\in\G_{(D)}$ with $c\ne0$. Then 
\be\label{Cgamint} C_\g(x)    \= \begin{cases} {\displaystyle H(-d/c) \m 
 |c|\,\v(\g)\,\sum_{k=1}^\infty\ch(k)\,\int_{k\,\g x}^\infty \frac{S(t)\,dt}{(ct-ak)^2}}
& \text{ if $x<-d/c$,} \\ {\displaystyle H(-d/c) \m
 |c|\,\v(\g)\,\sum_{k=1}^\infty\ch(k)\,\int_{-\infty}^{k\,\g x} \frac{S(t)\,dt}{(ak-ct)^2}}
& \text{ if $x>-d/c$\,.} \end{cases} \ee
 \end{proposition}

To see that this implies the continuity, we note that we can rewrite the right-hand side by
interchanging the summation and integration in the form (generalizing~\eqref{intC})
\be C_\g(x)  \=  H\bigl(-d/c\bigr) \m |c|\,\v(\g)\,\int_{-\infty}^\infty\,
  S(t)\,\Biggl(\sum_{k>0,\, k\g x\,\gtrless\,t} \frac{\ch(k)}{(ak-ct)^2}\Biggr) \,dt\ee
for $x\gtrless-d/c$.  The right-hand side of this formula
is continuous because $S(t)$ is piecewise continuous and bounded and has 
support in $|t|\ge1$, while the second factor of the integrand is piecewise continuous 
and bounded by a constant times~$t^{-2}$, as one sees easily by partial summation. Another
argument is that \eqref{Cgamint} is equivalent to the equality of distributions
\be\label{Cgamderiv} C'_\g(x) 
  \= -\,\v(\g)\,|c|\,\sgn(x+d/c)\sum_{k=1}^\infty\ch(k)\,\frac{S(k\,\g x)}{k}
  \= -\,\v(\g)\,c\;\sgn(cx+d)\,H\Bigl(\frac{ax+b}{cx+d}\Bigr)\;, \ee
(generalizing ~\eqref{Cderiv}) and we already know that the distribution~$H$ is locally the derivative of a continuous function.

For the proof of \eqref{Cgamint} we will restrict to the case when the matrix~$\g$ belongs 
to~$\G(D)$, which is sufficient for the proof of Theorem~\ref{QMFforD} because $\g\mapsto C_\g$ is a cocycle 
and we already know the continuity of $C_\g$~($=C$) for the special case $\g=S$ representing the non-trivial 
coset of $\G(D)$ in~$\G_{(D)}$.  To avoid distracting case distinctions, we consider only the case $x+d/c<0$ 
and $\g x>0$. The argument is the same in principle in all other cases, but we are actually working with 
ordered tuples of points on the 1-manifold $\mathbb P^1(\R)$, on which the group $\G_{(D)}$ acts in an 
orientation-preserving way, and it is notationally simpler to fix the positions of the points occurring with 
respect to the point at infinity, so that we can work on~$\R$ instead.  We may also assume (by replacing $\g$ 
by $-\g$ if necessary, but then remembering that $\g$ may be congruent to $-\bold1_2$ rather than $\bold1_2$ modulo~$D$)
that $c>0$. Finally, to be able to work with absolutely convergent sums and to reorder the terms freely, we 
will replace $H(x)$ by the function $H_s(x)$ ($s>1$) introduced in Section~\ref{Sec4.2} and $C_\g$ by the 
corresponding cocycle
\be\label{defCgams} C_{\g,s}(x) \=  H_s(x) \m \v(\g)\,|cx+d|^s\,H_s\Bigl(\frac{ax+b}{cx+d}\Bigr)\,, \ee
and only set $s=1$ at the end. The equation we have to prove then becomes
\be\label{Cgamints} C_{\g,s}(x)    \= H_s(-d/c) \m 
 c\,s\,\sum_{k=1}^\infty\ch(k)\,\int_{k\,\g x}^\infty \frac{S(t)\,dt}{(ct-ak)^{s+1}}\;. \ee
We denote the last term (without the minus sign) by~$A$.  Since the variable $t$ in the integral is always 
positive (because we are considering the case~$\g x>0$), we can replace $S(t)$ by its definition~\eqref{defS}
and interchange the summation and integration to getf
\bas A &\= \sum_{k,j>0}\ch(j)\,\ch(k)\,\int_{\max(j,k\,\g x)}^\infty \frac{c\,s\,dt}{(ct-ak)^{s+1}}  
  \=  \sum_{k,j>0} \frac{\ch(j)\,\ch(k)}{\max(cj-ak,k/|cx+d|)^s} \\
&\= \sums_{0<j\le k\g x}\frac{\ch(j)\,\ch(k)}{(k/|cx+d|)^s}
  \;+ \sums_{k>0,\;j\ge k\g x}\frac{\ch(j)\,\ch(k)}{(cj-ak)^s} \;.  \eas
In the second sum we replace the vector $\begin{pmatrix} j\\ k \end{pmatrix}$ by the vector
$-\g\begin{pmatrix} j\\ k \end{pmatrix}=\begin{pmatrix} -aj-bk\\ -cj-dk \end{pmatrix}$.
This does not change the product $\ch(j)\ch(k)$ because $\g$ is congruent to plus or minus the
identity modulo the period~$|D|$ of~$\ch$, changes the expression $cj-ak$ in the denominator to~$k$,
and changes the inequalities $k>0,\,j\ge k\g x$ to $kx\le j<-kd/c$ (which imply $k>0$).  Hence
\bas A &\= |cx+d|^s\,\sum_{k=1}^\infty \frac{\ch(k)}{k^s}\,S(k\g x)
  \;+ \sum_{k=1}^\infty \frac{\ch(k)}{k^s}\,\bigl(S(-kd/c)\m S(kx)\bigr) \\
   &\= H_s(-d/c) \m C_{\g,s}(x)\,  \eas
completing the proof of equation~\eqref{Cgamints} and hence of the proposition and theorem.

\medskip


\section{Modular forms are everywhere} \label{sec:MF}
As one can already see from the examples in the original article~\cite{QMF} where this notion was introduced,
quantum modular forms are sometimes related to actual modular forms.  These modular forms may be either 
holomorphic or Maass forms, with the quantum modular form in the latter case being related to the ``periods" 
of Maass forms in the sense developed in~\cite{LZ} and~\cite{BLZ}.  The quantum modular form~$H(x)$ that we 
have been studying in this paper turns out to be of this latter type, with the associated Maass form being an 
Eisenstein series with eigenvalue~1/4 for the hyperbolic Laplace operator.  In this final section we explain 
how this works, first for the case $D=-4$ studied in the first four sections of this paper.  In that case the 
relevant modular group $\G_\th$ is generated by a translation and an inversion, so that for both the period 
theory of the Maass form~$u$ and the quantum modularity of~$H$ one needs only the functional equation 
of the associated $L$-series. Then at the end we indicate how the quantum modularity of $H$ for 
general~$D$ follows from the full theory of periods of Maass forms.

We start with the case $D=-4$, so that $\chi$ is the character given by~\eqref{chi4L4}.  The relevant
modular form here is the Maass Eisenstein series \def\EE{u}
\be\label{uDef}  u(z) \= y^{1/2}\sum_{n=1}^\infty \chi(n)\,d(n)\,
      K_0\Bigl(\frac{\pi ny}2\Bigr)\,\sin\Bigl(\frac{\pi nx}2\Bigr) \qquad (z=x+iy\in\HH)\,, \ee
where $\HH$ denotes the upper half-plane and $K_0$ the usual $K$-Bessel function of order~0.  This is an
eigenfunction with eigenvalue 1/4 with respect to the hyperbolic Laplace operator 
$\D=-y^2\,\bigl(\frac{\p^2}{\p x^2}\+\frac{\p^2}{\p y^2}\bigr)$ and is a
modular function with character $\v$ for the group~$\G_\th$ defined in Section~\ref{sec.intro}, meaning
that $u(\g z)=\v(\g)u(z)$ for all $\g\in\G_\th$ or, more explicitly, that
\be\label{Einvariance} u(z+2) \= -\,u(z) \= u(-1/z) \qquad(z\in\HH). \ee
To see this, we observe that $u(z)$ is proportional to $E\bigl(\frac{z+1}4,\h\bigr)-E\bigl(\frac{z-1}4,\h\bigr)$,
where $E(z,s)$ is the usual non-holomorphic Eisenstein series of weight~0 and eigenvalue $s(1-s)$ for the
Laplace operator with respect to the full modular group~$\G_1=\SL$. (This follows easily from the well-known Fourier
expansion of $E(z,\h)$ as a linear combination of the three functions $\sqrt y$, $\sqrt y\log y$, and 
$\sqrt y\sum_{n\ne0} d(n)K_0(2\pi|n|y)e^{2\pi inx}$.) The transformation equations~\eqref{Einvariance} then 
follow from the $\SL$-invariance of $E(z,\h)$, the first one trivially since $E(z+1,s)=E(z,s)$ and the
second by using the invariance of $E(z,s)$ under $\sma 10{\mp4}1$ to get 
$$ E\Bigl(\frac{-1/z\pm1}4,\,\frac12\Bigr) \= E\Bigl(\frac{(-1/z\pm1)/4}{\mp4((-1/z\pm1)/4)+1},\,\frac12\Bigr)
 \= E\Bigl(\frac{z\mp1}4,\,\frac12\Bigr)\;.  $$
 
We now associate to~$u(z)$ the periodic holomorphic function $f$ on $\CC\sm\R=\HH^+\cup\HH^-$ 
(where $\HH^\pm=\{z\in\CC\mid\pm\Im(z)>0\}$) having the same Fourier coefficients as~$u$, i.e.,
\be\label{fDef} f(z) \= \sum_{n=1}^\infty \ch(n)\,d(n)\,q^{\pm n/4} \qquad(z\in\HH^\pm,\; q=e^{2\pi iz}\bigr)\,. \ee
\begin{proposition} \label{prop:psi4} The  {\it period function} $\psi(z)$ defined by
\be\label{psiDef} \psi(z)\= f(z) \+ \frac1z\,f\Bigl(-\frac1z\Bigr)\qquad(z\in\CC\sm\R) \ee
extends holomorphically from $\CC\sm\R$ to $\CC'=\CC\sm(-\infty,0]\,$.
 \end{proposition}

\noindent{\it Proof.} We follow the proof of the corresponding result given in Chapter~1 of~\cite{LZ} for 
Maass cusp forms on the full modular group~$\G_1$. (See Theorem on p.~202 of~\cite{LZ}, which also gives
a converse statement characterizing cusp forms in terms of holomorphic functions $\psi$ on~$\CC'$ satisfying
a certain functional equation.)  That proof required only the functional equation of the $L$-series associated
to the Maass cusp form, which worked because the group $\G_1$ is generated by the translation~$T$ and the
inversion~$S$, and can be applied here because  $\G_\th=\la S,T^2\ra$ has a similar structure. Our situation
is also a little different because our function $u(z)$ is an Eisenstein series rather than a cusp form, but
since the Fourier expansion~\eqref{uDef} of~$u(z)$ at infinity has no constant term this has no effect on the proof.

The argument in~\cite{LZ} was first to write the $L$-series of the Maass form $u(z)$ as the Mellin transform 
of the restriction of~$u$ (or of its normal derivative in the case of an odd cusp form) to the imaginary
axis, multiplied by a suitable gamma factor, and to deduce from this relationship and the $S$-invariance of~$u$
a functional equation for the $L$-series.  One then observed that the Mellin transforms of the restrictions to 
the positive or negative imaginary axis of both the associated periodic holomorphic function~$f$ and the associated
period function~$\psi$ were also equal, up to different gamma factors, to the same $L$-series, and the
functional equation of this $L$-series combined with Mellin inversion then led to a formula for~$\psi$
that applied in all~$\CC'$ rather than just on~$\CC\sm\R$.  Here the first step, which would require the normal
derivative since the restriction of the function~\eqref{uDef} to~$i\R_+$ vanishes, can be skipped since
the $L$-series $L(u;s)=\sum \ch(n)d(n)n^{-s}$ is simply the square of the Dirichlet $L$-series~\eqref{chi4L4}
and hence has a known functional equation.  The Mellin transforms of the restrictions to $f$ and $\psi$
to the positive or negative imaginary axis are then given by
$$ \wt f_\pm(s) \: \int_0^\infty f(\pm iy)\,y^{s-1}\,dy \= \frac{\G(s)}{(\pi/2)^s}\,L(s)^2\qquad(\Re(s)>0),$$
\bas \wt \psi_\pm(s) &\= \int_0^\infty \psi(\pm iy)\,y^{s-1}\,dy 
\= \int_0^\infty \Bigl(f(\pm iy) \,\mp\, \frac iy\,f(\pm i/y)\Bigr)\,y^{s-1}\,dy \\
&\= \wt f_\pm(s) \,\mp\,i\,\wt f_\pm(1-s) \= \frac{\G(s)\,L(s)^2}{(\pi/2)^s\,\cos(\pi s/2)}\,e^{\mp i\pi s}
 \qquad(0<\Re(s)<1), \eas
where in the last line we have used the functional equation of~$L(s)$ and standard identities for gamma
functions. By the Mellin inversion formula we deduce that
$$ \psi(\pm iy)\=\frac1{2\pi i}\,\int_{\Re(s)=c} \frac{\G(s)\,L(s)^2}{(\pi/2)^s\,\cos(\pi s/2)}\,(\pm iy)^{-s}\,ds 
\qquad(y>0) $$
for any~$c\in(0,1)$, and analytic continuation from $i\R\sm\{0\}$ to $\CC\sm\R$ then gives the formula
$$ \psi(z) \= \frac1{2\pi i}\,\int_{\Re(s)=c} \Bigl[\frac{\G(s)\,L(s)^2}{(\pi/2)^s\,\cos(\pi s/2)}\Bigr]\,z^{-s}\,ds $$
for all $z\in\CC\sm\R\,$.  But the factor in square brackets is bounded by a power of~$s$ times $e^{-\pi|s|}$
for $|s|\to\infty$ on the vertical line $c+i\R$, so the integral on the right-hand side of this equation is 
absolutely convergent for all $z\in\CC$ with $|\arg(z)|<\pi$, i.e., for all~$z\in\CC'\,$.   \endpf

We now show how to connect the functions $f(z)$ and $\psi(z)$ to the functions $H(x)$ and~$C(x)$, respectively,
and to deduce the continuity of~$C$ -- and hence the quantum modularity of~$H$ -- from Proposition~\ref{prop:psi4}
for the period function~$\psi$. The first step is easy, since from the Fourier expansions~\eqref{fDef} 
and~\eqref{HFour} and standard integral formulae we find the relation
$$  f(z) \;\doteq\; \int_{-\infty}^\infty \frac{H(t)}{(t-z)^2}\,dt \qquad(z\in\CC\sm\R) $$
between the periodic function $f(z)$ and the periodic distribution $H(x)$. (Here and in what follows the symbol 
$\,\doteq\,$ denotes equality up to easily computed scalar factors whose
values are irrelevant for the argument.) Replacing $t$ by $1/t$ in the integral, we obtain
\be\label{fint}  f(z) \;\doteq\; \int_{-\infty}^\infty \frac{C'(t)}{(1-tz)^2}\,dt 
\;\doteq\; z\,\int_{-\infty}^\infty \frac{C(t)}{(1-tz)^3}\,dt \,, \ee
where we have used~\eqref{Cderiv} to get the first equality and integration by parts for the second.
But now replacing $t$ by $1/t$ and using the functional equation $C(t)=|t|C(1/t)$ we find
\be  \frac1z\,f\Bigl(-\frac1z\Bigr) \= \frac1z\,f\Bigl(\frac1z\Bigr) 
\;\doteq\; -\,z\,\int_{-\infty}^\infty \frac{\sgn(t)\,C(t)}{(1-tz)^3}\,dt \ee
with the same proportionality constant as in~\eqref{fint}, and adding these two equations gives the relation
\be\label{psifromC} \psi(z) \;\doteq\;  z\,\int_{-\infty}^0 \frac{C(t)}{(1-tz)^3}\,dt 
\;\doteq\;  z\,\int_{-\infty}^0 \frac{C(t)}{(t-z)^3}\,dt   \qquad(z\in\CC\sm\R) \ee
between the period function of the Maass wave form~$u$ and the function~$C(x)$.  This establishes
the desired connection between the continuity of $C(x)=H(x)+|x|H(1/x)$, which expresses the quantum modularity
of~$H$, and Proposition~\ref{prop:psi4}, which expresses the modularity of~$u$: in one direction,
if we know that $C(x)$ is continuous (and bounded by $\,\min(1,|x|)$), then~\eqref{psifromC} immediately
gives the analytic continuation of~$\psi(z)$ from~$\CC\sm\R$ to~$\CC'$; and conversely, the
inversion formula for the Stieltjes transform given in~\cite{Sch} lets us invert~\eqref{psifromC} to get
\be\label{Cfrompsi} C(t) \;\doteq\; t\, \int_{\mathcal C} \, \psi(tz) \frac{1 + z}z\,dz \ee
where ${\mathcal C}$ is any contour with endpoints at $z = -1$ which encloses the origin, so that the
continuity of~$C$ is a direct consequence of the holomorphy of~$\psi$ in~$\CC'$.

We observe that the entire discussion given here applies in an essentially unchanged form to the more general
functions $C_s$ and $H_s$ discussed in~\S\ref{Sec4.2}, with the Maass form~\eqref{uDef} replaced by the form
$E\bigl(\frac{z+1}4,\frac s2\bigr)-E\bigl(\frac{z-1}4,\frac s2\bigr)$ with spectral parameter~$s/2$ instead of~$1/2$.

Finally, we consider the case when the character~$\ch$ in~\eqref{chi4L4} is replaced by an arbitrary
primitive odd Dirichlet character $\ch_D$.  Then all of the above calculations still go through: the
function~$u$ is defined by equation~\eqref{uDef} with $\ch$ replaced by~$\ch_D$ and $\pi/2$
by $2\pi/|D|$, which is again a Maass form (with spectral parameter~$\h$ and character~$\v$ on
the group~$\G_{(D)}$) because it is proportional to $\sum_{r\!\!\pmod D}\ch(r)E\bigl(\frac{z+r}{|D|},\h\bigr)$;
the associated periodic function~$f$ and period function~$\psi$ are defined by~\eqref{fDef} (with the new~$\ch$
and with $q^{n/4}$ replaced by~$q^{n/|D|}$) and~\eqref{psiDef} (with no change at all); and 
Proposition~\ref{prop:psi4} remains true with the same proof.  The difference, however, is that this
proposition is no longer equivalent to the modularity of~$u$, but only to its invariance (up to sign) under the
transformations $S$ and $T^D$ (or $T^{D/2}$ if $D$ is even), which in general generate a subgroup of~$\G_{(D)}$
of infinite order, as already discussed in Section~\ref{sec:OtherChi}.  To get the full modularity (of this~$u$
or any other potential Maass form on~$\G_{(D)}$), we need to generalize Proposition~\ref{prop:psi4} to the
statement that for any matrix $\g=\sma abcd\in\G_{(D)}$ the function
$$ \psi_\g(z)\= f(z) \m  \frac{\v(\g)}{cz+d}\,f\Bigl(\frac{az+b}{cz+d}\Bigr)  $$
extends holomorphically from $\CC\sm\R$ to $\CC\sm(-\infty,-d/c]\,$ if $c>0$ or $\CC\sm[-d/c,\infty)$ if $c<0$.
This statement can be proved in several ways and can be linked by a discussion similar to the one above to 
the continuity property of $C_\g(x)$ stated in Theorem~\ref{QMFforD}. But the simplest approach is to relate
the functions $C_\g(x)$ directly to the invariant distribution associated to~$u$ in the sense developed
in~\cite{L}, \cite{LZ} and~\cite{BLZ}. Specifically, this theory says that an eigenfunction~$u$ of~$\D$ with spectral 
parameter~$s$ and Fourier expansion $u(x+iy)=\sqrt y\,\sum_{n\ne0}A_n\,K_{s-1/2}(\l|n|y)e^{i\l nx}$ is invariant 
(possibly with character) under the action of a Fuchsian group~$\G$ if and only if the associated distribution 
$U(x)=\sum_{n\ne0}|n|^{s-1/2}\,A_n\,e^{i\l nx}$ on~$\mathbb P^1(\R)$ is invariant (with the same character) with respect 
to the group action $(U|\g)(x)=|cx+d|^{-2s}U(\g x)$ \hbox{for~$\g=\sma abcd\in\G$}. (Here the word ``distribution" must 
be interpreted correctly, namely, as a functional on the space of test functions $\phi(x)$ that are smooth on~$\R$ and 
for which $|x|^{-2s}\phi(1/x)$ is smooth near~$x=0$.) In our case the distribution~$U$ associated to~$u$ is given by
$U(x)=\sum_{n=1}^\infty \ch(n)d(n)\sin(2\pi nx/|D|)$ and is related to the distribution~\eqref{newHFour} by 
$U(x)\doteq H'(x)$.  But that means that if we differentiate the definition~\eqref{defCgam} of~$C_\g(x)$ to get
$$  C_\g'(x) \= H'(x) \m \frac{\v(\g)}{|cx+d|}\,H'(\g x) \m \v(\g)\,|c|\,\sgn(x+d/c)\,H\g x)\,, $$
then the first two terms on the right cancel and we recover equation~\eqref{Cgamderiv}.



\enlargethispage{3mm}

\end{document}